\documentclass[11pt, reqno]{amsart}

\usepackage[T1]{fontenc}

\usepackage[utf8]{inputenc}

\usepackage[margin=1in]{geometry}
\usepackage{amsthm,amsmath,amssymb}

\usepackage{tikz,xcolor,graphicx}
\usepackage{tikz-cd}

\usepackage{extarrows}

\usepackage{fancyhdr}

\usepackage{mathtools}
\usepackage{thmtools}
\usepackage{thm-restate}

\usepackage{booktabs}

\usepackage[noadjust]{cite}

\usepackage{listings}
\usepackage[shortlabels]{enumitem}

\usepackage{comment}

\usepackage{nicefrac}

\usepackage[unicode]{hyperref}
\hypersetup{colorlinks=true,urlcolor=black,citecolor=blue,linkcolor=blue}
\usepackage[noabbrev,capitalize]{cleveref}

\usepackage{url}

\usepackage[color, final]{showkeys} 

\colorlet{refkey}{orange!20}
\colorlet{labelkey}{blue!80}

\usepackage[colorinlistoftodos]{todonotes}

\usepackage{marginnote}

\usetikzlibrary{calc}

\usetikzlibrary{arrows, shapes} 
\tikzstyle{vertex}=[circle,fill=black!25,minimum size=20pt,inner sep=0pt]
\tikzstyle{selected vertex} = [vertex, fill=red!24]
\tikzstyle{edge} = [draw,thick,-]
\tikzstyle{weight} = [font=\small]
\tikzstyle{selected edge} = [draw,line width=5pt,-,red!50]
\tikzstyle{ignored edge} = [draw,line width=5pt,-,black!20]

\usepackage{float}
\usepackage{pgfplots}
\usepgfplotslibrary{fillbetween}
\graphicspath{ {images/} }
\pgfplotsset{compat=1.16}

\theoremstyle{plain}
\newtheorem{theorem}{Theorem}[section]

\newtheorem{proposition}[theorem]{Proposition}
\newtheorem{lemma}[theorem]{Lemma}

\newtheorem{corollary}[theorem]{Corollary}
\newtheorem{conjecture}[theorem]{Conjecture}

\newtheorem{assumption}[theorem]{Assumption}

\newtheorem{innercustomcor}{Corollary}
\newtheorem*{question}{Question}
\newtheorem*{corollary-nonum}{Corollary}

\theoremstyle{definition}
\newtheorem{definition}[theorem]{Definition}

\theoremstyle{remark}
\newtheorem*{remark}{Remark}

\newenvironment{customcor}[1]
  {\renewcommand\theinnercustomcor{#1}\innercustomcor}
  {\endinnercustomcor}

\usepackage{titletoc}
\newif\ifnotes
\notestrue

\newcommand{\NN}{\mathbb{N}}

\newcommand{\ZZ}{\mathbb{Z}}

\newcommand{\Fp}{\mathbb{F}_p}

\newcommand{\Q}{\mathbb{Q}}

\newcommand{\R}{\mathbb{R}}
\renewcommand{\C}{\mathbb{C}}
\newcommand{\Z}{\mathbb{Z}}

\newcommand{\lcm}{\operatorname{lcm}}

\renewcommand{\O}{\mathcal{O}}

\newcommand{\cO}{\mathcal{O}}

\newcommand{\cR}{\mathcal{R}}

\renewcommand{\t}{\text}
\newcommand{\f}{\frac}
\newcommand{\eps}{\epsilon}
\renewcommand{\Re}{\mathrm{Re}}
\renewcommand{\Im}{\mathrm{Im}}

\newcommand{\textdef}[1]{\emph{#1}}
\newcommand{\argument}{\,{\boldsymbol\cdot}\,}

\newcommand{\paren}[1]{\left( #1 \right)}

\newcommand{\abs}[1]{\left\lvert #1 \right\rvert}

\newcommand{\absBigg}[1]{\Bigg\lvert #1 \Bigg\rvert}

\newcommand{\ceil}[1]{\left\lceil #1 \right\rceil}

\newcommand{\interior}[1]{%
  {\kern0pt#1}^{\mathrm{o}}%
}

\newcommand\restr[2]{{
  \left.\kern-\nulldelimiterspace 
  #1 
  \vphantom{\big|} 
  \right|_{#2} 
}}

\makeatletter
\newcommand{\customlabel}[2]{%
   \protected@write \@auxout {}{\string \newlabel {#1}{{#2}{\thepage}{#2}{#1}{}} }%
   \hypertarget{#1}{}
}
\makeatother

\newcommand{\restatableTheorem}[4]{
    \begin{#2}
        \label{#3}
        #4
    \end{#2}
    \theoremstyle{plain}
    \newtheorem*{{restatedtheorem}#3}{\cref*{#3}}
    \newcommand{#1}{
        \begin{{restatedtheorem}#3}
            \customlabel{#3:restated}{\labelcref*{#3}}
            #4
        \end{{restatedtheorem}#3}
    }
}

\newcommand{\crefrestated}[1]{\namecref{#1}~\ref{#1:restated}}

\DeclareMathOperator{\sqf}{\mathrm{sqf}}
\DeclareMathOperator{\sgn}{\mathrm{sgn}}
\DeclareMathOperator{\rad}{\mathrm{rad}}   
\renewcommand{\mod}[1]{\ (\mathrm{mod}\ #1)}
\newcommand{\leg}[2]{\paren{\tfrac{#1}{#2}}}
\newcommand{\nf}[2]{\nicefrac{#1}{#2}}
\newcommand{\tf}[2]{\tfrac{#1}{#2}}
\newcommand{\veps}{\varepsilon}
\usepackage{etoolbox}

\makeatletter
\patchcmd{\@settitle}{\uppercasenonmath\@title}{\Large}{}{}
\patchcmd{\@setauthors}{\MakeUppercase}{\large}{}{}
\makeatother

\title[Quadratic Fields Admitting Elliptic Curves with Good Reduction Everywhere]{Quadratic Fields Admitting Elliptic Curves with\\ Rational \texorpdfstring{$j$}{j}-Invariant and Good Reduction Everywhere}

\author{Benjamin Matschke}
\address{~}
\email{matschke@bu.edu}
\urladdr{\url{https://math.bu.edu/people/matschke/}}
\thanks{The first author was supported by Simons Foundation grant \#550023.}

\author{Abhijit S. Mudigonda}
\address{~}
\email{abhijitm@uchicago.edu}
\urladdr{\url{https://abhijit-mudigonda.github.io/math/}} 

\begin{document}
\begin{abstract}
Clemm and Trebat-Leder (2014) proved that the number of quadratic number fields with absolute discriminant bounded by $x$ over which there exist elliptic curves with good reduction everywhere and rational $j$-invariant is $\gg x\log^{\nicefrac{-1}{2}}(x)$.
In this paper, we assume the $abc$-conjecture to show the sharp asymptotic $\sim cx\log^{\nicefrac{-1}{2}}(x)$ for this number, obtaining formulae for $c$ in both the real and imaginary cases. Our method has three ingredients:

\begin{enumerate}
    \item We make progress towards a conjecture of Granville: 
Given a fixed elliptic curve $E/\Q$ with short Weierstrass equation $y^2 = f(x)$ for reducible $f \in \Z[x]$, we show that the number of integers~$d$, $|d| \leq D$, for which the quadratic twist $dy^2 = f(x)$ has an integral non-$2$-torsion point is at most $D^{\nf{2}{3}+o(1)}$, assuming the $abc$-conjecture. 
    \item We apply the Selberg--Delange method to obtain a Tauberian theorem which allows us to count integers satisfying certain congruences while also being divisible only by certain primes. 
    \item We show that for a polynomially sparse subset of the natural numbers, the number of pairs of elements with least common multiple at most $x$ is $O(x^{1-\eps})$ for some $\eps > 0$. We also exhibit a matching lower bound. 
\end{enumerate}
If instead of the $abc$-conjecture we assume a particular tail bound, we can prove all the aforementioned results and that the coefficient $c$ above is greater in the real quadratic case than in the imaginary quadratic case, in agreement with an experimentally observed bias. 
\end{abstract}

\maketitle 
\vspace*{-\baselineskip}

\section{Introduction}
\label{sec:intro}

An elliptic curve defined over a number field $K$ is said to have \textdef{good reduction everywhere} if it has good reduction at every prime ideal of the ring of integers of $K$. Tate showed that there are no such elliptic curves over $\Q$~\cite{ogg}, but this is not the case for all fields. For example, over $K := \Q(\sqrt{29})$ the elliptic curve 
\begin{equation*}
    y^2 +xy + \paren{\f{5+\sqrt{29}}{2}}^2y = x^3
\end{equation*}
has good reduction everywhere. The existence and properties of quadratic fields admitting such elliptic curves has been studied extensively~\cite{comalada-fixedj, comalada-nart, comalada-quad, kagawa-cubdisc, kida-potential, kida-nonexistence, stroeker-egr, setzer-2, zhao, kida-certain, kida-imaginary, ishii-nonex}, and extensions to higher-degree fields were considered by Takeshi~\cite{takeshi-cub, takeshi-gen}. Algorithms for computing such elliptic curves were given by Kida~\cite{kida-comp} for quadratic fields and by Cremona and Lingham~\cite{cremona-lingham}, Koutsianas~\cite{koutsianas}, and the first author~\cite{matschke} over general number fields. 
We will be interested in the following statistical question:

\begin{question}
    How often does a real (resp.\ imaginary) quadratic field admit an elliptic curve with good reduction everywhere and rational $j$-invariant?
\end{question}

As with any such question, we must define ``often'' with respect to an ordering of quadratic fields. Let $R(x)$ (resp.\ $I(x)$) be the number of real (resp.\ imaginary) quadratic fields $K/\Q$ with discriminants of absolute value at most $x$ and over which there exist elliptic curves with rational $j$-invariant and with good reduction at every prime of $K$. Setzer~\cite{setzer-crit} gave an explicit criterion for $m$ such that $\Q(\sqrt{m})$ admits an elliptic curve with rational $j$-invariant and good reduction everywhere. Coupling this criterion with a lower bound of Serre~\cite{serre} on the sizes of particular sifted sets of integers, Clemm and Trebat-Leder~\cite{ctl} gave a lower bound on the quantities of interest. 
\begin{theorem}[Clemm and Trebat-Leder~\cite{ctl}]
    \label{thm:ctl}
    \begin{equation*}
        R(x) \gg \f{x}{\sqrt{\log x}} \quad\t{and}\quad I(x) \gg \f{x}{\sqrt{\log x}}.
    \end{equation*}
\end{theorem}

\subsection{Main results}

Our main theorem is a sharp result for the asymptotic behavior of $R(x)$ and $I(x)$ assuming the $abc$-conjecture. Because the expressions for the constants are somewhat technical, we defer some of their descriptions to later in the paper. 

\restatableTheorem{\rixconst}{maintheorem}{thm:rix-const}{
    Assuming the $abc$-conjecture,
    \[R(x) \sim \f{c_Rx}{\sqrt{\log x}} \quad\t{and}\quad  I(x) \sim \f{c_{I}x}{\sqrt{\log x}}\]
    where 
    \[c_R = \sum_{\substack{d \in \Z \\ d \t{ good}}} \f{c_dc'_{d,R}}{|d|2^{\omega(d)}} \quad\t{and}\quad c_I = \sum_{\substack{d \in \Z \\ d \t{ good}}} \f{c_dc'_{d,I}}{|d|2^{\omega(d)}},\]
    where the set of good $d$ is defined in~\eqref{eq:good}, $c_d$ is as in~\cref{cor:ridx-asymp} and $c'_{d,R}$ and $c'_{d,I}$ are as in~\cref{lem:ridx-const}.
}

We may worry that the condition of ``rational $j$-invariant'' is too restrictive. However, conditionally complete\footnote{The tables are complete assuming the generalized Riemann hypothesis.} tables of elliptic curves with bounded absolute discriminant and good reduction everywhere due to the first author~\cite{matschke} suggest that most quadratic fields admitting any curves with good reduction everywhere admit at least one such curve with rational $j$-invariant. We discuss the possibility of extending these results to all $j$-invariants further in~\cref{subsec:future}.

One ingredient in the proof of \cref{thm:rix-const} is an upper bound on how often quadratic twists of some elliptic curves over $\Q$ have integral points. Consider an elliptic curve $E/\Q$ with short Weierstrass equation $y^2 = f(x)$. We are interested in counting the number of quadratic twists $E_d$ of $E$ -- with models $dy^2 = f(x)$ -- that have integral points. Note that if $E$ has a two-torsion point, say $(a:0:1)$, then every twist $E_d$ will also have the two-torsion point $(a:0:1)$. As such, we restrict our attention to determining the existence of a \textdef{nontrivial} integral point on $E_d$ -- an integral point that is not two-torsion.

Formally, we prove an upper bound on the quantity $|T_E(D)|$, where
\begin{equation}
    \label{eq:ted}
    T_E(D) := \{d \in \Z \colon |d| \leq D, d \t{ is squarefree and } E_d \t{ has a nontrivial integral point}\}.
\end{equation}

Granville~\cite{granville-twists} showed the following conditional upper bound on the analogous quantity for hyperelliptic curves.  
\begin{theorem}[Granville~\cite{granville-twists}]
    \label{thm:granville-1}
    Assume that the $abc$-conjecture is true. Let $C$ be a hyperelliptic curve given by the integral model $y^2 = f(x)$ where $f \in \Z[x]$ has degree at least three (i.e.\ the genus is at least one) and is separable. Then, 
    \begin{equation*}
        |T_C(D)| \leq D^{\f{1}{\deg f-2} + o(1)}.
    \end{equation*}
\end{theorem}

We will be interested in the case where $\deg f = 3$, in which case~\cref{thm:granville-1} is trivial. 
However, Granville also conjectures that the $\deg f - 2$ in the denominator of the exponent can be replaced with~$\deg f$.

\begin{conjecture}[Granville~\cite{granville-twists}]
    \label{conj:granville-1}
    Let $C$ be a hyperelliptic curve given by the integral model $y^2 = f(x)$ where $f \in \Z[x]$ has degree at least three (i.e.\ the genus is at least one) and is separable. Then, 
    \begin{equation*}
        |T_C(D)| \sim \kappa_fD^{\f{1}{\deg f} + o(1)},
    \end{equation*}
    where $\kappa_f$ is some constant that can be determined explicitly given $f$.
\end{conjecture}

Granville proves~\cref{conj:granville-1} for polynomials $f$ of degree at least $7$ that split into linear factors over~$\Q$. We make progress towards~\cref{conj:granville-1} when the degree of $f$ is $3$.

\restatableTheorem{\redec}{maintheorem}{thm:red-ec}{
    Assume that the $abc$-conjecture is true. Let $E$ be an elliptic curve over $\Q$ with short Weierstrass equation $y^2 = f(x)$ where $f(x) \in \Z[x]$ is reducible\footnote{Equivalently, $E$ has a rational Weierstrass point.} over~$\Q$. Then,
    \begin{equation*}
        |T_E(D)| \leq D^{\nf{2}{3} + o(1)}. 
    \end{equation*}
}

We may apply~\cref{thm:red-ec} to the set of good $d$ in~\cref{thm:rix-const} to show (among other things) that the sum in the statement of~\cref{thm:rix-const} converges. At a high level, this sum is actually a union bound over the good $d \in \Z$, and thus it yields an upper bound on the leading constant in~\cref{thm:rix-const}. To show the matching lower bound, we control the second term in the inclusion-exclusion sequence via a result on the sizes of pairwise least common multiples of ``polynomially sparse'' subsets. This result may be of independent interest. 

\restatableTheorem{\lcmsparse}{maintheorem}{thm:lcm-sparse}{
    A set $S \subseteq \NN$ of squarefree numbers is called $\beta$-polynomially sparse if for $\beta \in (0,1)$ we have
    \[\#\{n \leq x \colon n \in S\} \leq x^{1-\beta+o(1)}\]
    as $x$ approaches $+\infty$. For any such $S$, the set
    \[\{(n, n') \colon n \in S, n' \in S, \lcm(n,n') \leq x\}.\]
    is $\f{\beta}{2-\beta}$-polynomially sparse. Furthermore, there are sets for which this is tight. 
}

\begin{corollary-nonum}
    The set of pairwise least common multiples of a polynomially sparse set of squarefree numbers is polynomially sparse.
\end{corollary-nonum}

The initial motivation for studying the constants $c_R$ and $c_I$ in~\cref{thm:rix-const} was the observation that most quadratic fields admitting curves with everywhere good reduction appear to be real. Assuming the $abc$-conjecture, we are able to prove numerical lower bounds on $c_R$ and $c_I$. 

\restatableTheorem{\constlb}{maincorollary}{cor:const-lb}{
    Assuming the $abc$-conjecture,~\cref{thm:rix-const} holds for
    \[c_R \geq 0.1255 \quad\t{and}\quad c_I \geq 0.01109.\]
}

We expect these values, obtained by evaluating the sum in~\cref{thm:rix-const} for many good $d$, to be very close to the truth.
Indeed, the aforementioned tables of elliptic curves show that,
under the generalized Riemann hypothesis, $R(20000)=728$ and $I(20000)=97$.
\footnote{These numbers of fields become $852$ and $97$, respectively, if we drop the condition of rational $j$-invariants.}
Writing $\tilde{c}_R$ and $\tilde{c}_I$ to denote the constants in~\cref{cor:const-lb},
we have $\tilde{c}_R\, x_0\log^{\nicefrac{-1}{2}}x_0 \approx 797$ and $\tilde{c}_I\, x_0\log^{\nicefrac{-1}{2}}x_0 \approx 70$ for $x_0 = 20000$,
roughly in line with the true values. 

\begin{restatable}{maincorollary}{constub}
    \label{cor:const-ub}
    Let $E$ be the elliptic curve given by the short Weierstrass equation $y^2 = x^3-1728$. Assume, instead of the $abc$-conjecture, that $|T_E(D)| \leq 5 D^{0.35}$. Then,~\cref{thm:rix-const} holds with
    \[0.1255 \leq c_R \leq 0.1489 \quad\t{and}\quad 0.01109 \leq c_I \leq 0.03446.\]
    In particular, $c_R > c_I$ under this hypothesis.  
\end{restatable}

Experimentally, we have checked (\cref{fig:d-conj}) that this hypothesis holds comfortably for all $D \leq 10000$. We also motivate this hypothesis using the aforementioned conjecture of Granville (\cref{conj:granville-1}).

\subsection{Techniques and an overview of the proofs}
\label{subsec:tech}

The first input used in proving~\cref{thm:rix-const} is a criterion of Setzer~\cite{setzer-crit} for when $\Q(\sqrt{m})$ admits an elliptic curve with good reduction everywhere and rational $j$-invariant. We state the criterion formally as~\cref{thm:setzer-crit} but at a high level, it tells us that $\Q(\sqrt{m})$ admits such an elliptic curve if and only if $m$ can be factored as $nd$ for some $d$ such that 
\begin{enumerate}[(a)]
    \item $d$ is the squarefree part of $r^3 - 1728$ for $r$ in some positive density subset of the integers; \label{tech-a}
    \item $n$ is divisible only by primes satisfying certain quadratic residuosity conditions with respect to $d$. For any $d$, the set of primes satisfying this condition has natural density $\f{1}{2}$; \label{tech-b}
    \item the image of $n$ in $(\Z/4d\Z)^{\times}$ lies in a specified subset. \label{tech-c}
\end{enumerate}

We are interested in upper bounding $R(x)$, the number of real quadratic fields with discriminant at most $x$ and which satisfy conditions~\ref{tech-a},~\ref{tech-b}, and~\ref{tech-c} (the same techniques apply to $I(x)$).

We first show, assuming the $abc$-conjecture, that the set of $d \in \Z$ which satisfy~\ref{tech-a} is polynomially sparse. This is a corollary of~\cref{thm:red-ec}. The key idea motivating the proof of~\cref{thm:red-ec} is the relationship between squarefree parts and quadratic twists of elliptic curves. Consider~\ref{tech-a} above. We have that $d$ is the squarefree part of $r^3-1728$ if and only if for some integer $t$ we have $dt^2 = r^3-1728$. This happens if and only if the quadratic twist by $d$ of $E: t^2 = r^3-1728$ has an integral point. By definition, $T_E(D)$ (\eqref{eq:ted}) counts the nontrivial points and hence the number of nonzero $d$ which arise as squarefree parts of $r^3-1728$.

Next, we fix some $d$ satisfying~\ref{tech-a} and study the asymptotics of $R_d(x)$, the contribution to $R(x)$ from those $m$ which are divisible by this $d$. Then, we have
\begin{equation}
    \label{eq:tech-1}
    R(x) \leq \sum_{d \t{ sat.\,\ref{tech-a}}} R_d(x).
\end{equation}

This approach is motivated by the lower bound of Clemm and Trebat-Leder~\cite{ctl} (\cref{thm:ctl}). They chose a single value of $d$ satisfying~\ref{tech-a} and for which~\ref{tech-c} is trivial. A result of Serre~\cite{serre} implies a lower bound on the number of $n$ satisfying~\ref{tech-b} and shows that for this choice of $d$, $R_d(x) \gg \f{x}{\sqrt{\log x}}$. This may seem surprising, as it means that even without considering multiple values of $d$ they are already able to obtain the correct order of growth of $R(x)$! This happens because the set of $d$ satisfying~\ref{tech-a} is very sparse. Indeed, it turns out that for any such $d$, 

\begin{equation}
    \label{eq:ridx-asymp-inf}
    R_d(x) \ll \f{(1+o_d(x))c_dx}{\abs{d}\sqrt{\log x}},
\end{equation}
where $c_d$ grows very slowly as $|d|$ goes to infinity. 
We prove~\eqref{eq:ridx-asymp-inf} using Selberg-Delange theory (in particular, \cref{thm:selberg-delange}), which is also the general theory underlying the bound of Serre~\cite{serre} used in the work of Clemm and Trebat-Leder. Selberg-Delange theory gives us a Tauberian theorem for Dirichlet series which can be expressed as $\zeta(s)^{\rho}G(s)$ for $\rho \in \C$ and $G(s)$ holomorphic in a neighborhood around $s = 1$. For $\rho$ a nonzero real, it tells us that the sum of coefficients of the series up to $x$ is asymptotically $cx\log^{\rho-1} x$ for some explicit constant $c$ depending on $G$. To obtain our upper bound, we apply Selberg-Delange theory to the Dirichlet series $F(s)$ whose coefficients are the values of the characteristic function of~\ref{tech-b} -- since we just need an upper bound, it is fine to ignore~\ref{tech-c} for now. It turns out that because the number of ``valid'' primes in~\ref{tech-b} is half of all primes, $F(s) = \zeta(s)^{\nf{1}{2}}G(s)$ for some $G(s)$ holomorphic around $s = 1$. This then implies~\eqref{eq:ridx-asymp-inf}. Applying~\eqref{eq:ridx-asymp-inf} to~\eqref{eq:tech-1}, we deduce

\begin{equation}
    \label{eq:tech-15}
    R(x) \ll \sum_{d \t{ sat.\,\ref{tech-a}}} \f{(1+o_d(x))c_dx}{|d|\sqrt{\log x - \log |d|}}.
\end{equation}

Summation by parts shows that the sum of the reciprocals of the elements of a polynomially sparse set converges (\cref{lem:nat-to-harm}), and applying this to the set of $d$ satisfying~\ref{tech-a} (which is polynomially sparse by~\cref{thm:red-ec}), we have that the series
\[\sum_{d \t{ sat.\,\ref{tech-a}}} \f{c_d}{|d|}\]
converges. This allows us to uniformly bound the $o_d(x)$ terms in~\eqref{eq:tech-15} and obtain
\begin{equation}
    \label{eq:tech-175}
    R(x) \ll \f{x}{\sqrt{\log x}} \sum_{d \t{ sat.\,\ref{tech-a}}} \f{c_d}{|d|} \ll \f{x}{\sqrt{\log x}}.
\end{equation}

In order to compute the implicit constant in~\eqref{eq:tech-175}, we start by computing the implicit constant in~\eqref{eq:ridx-asymp-inf} (\cref{lem:ridx-const}). To do this, we fix a $d$ and count those $n$ which satisfy~\ref{tech-b} and~\ref{tech-c}. Selberg-Delange theory can be applied directly to obtain the exact constant if we are only interested in the Dirichlet series of~\ref{tech-b}. We need to study the Dirichlet series of \ref{tech-b}$\land$\ref{tech-c},\footnote{The Dirichlet series whose coefficient at $n$ is $1$ if and only if $n$ satisfies~\ref{tech-b} and~\ref{tech-c}} which is the (Rankin-Selberg) convolution of the Dirichlet series for~\ref{tech-b} and the Dirichlet series for~\ref{tech-c}. However,~\ref{tech-c} need not be a multiplicative property, and hence its Dirichlet series need not have an Euler product. We address this by expressing it as a linear combination of Dirichlet $L$-series and noting that only one term of the linear combination contributes to the overall asymptotics. We then apply~\cref{thm:selberg-delange} to this term to obtain our result on $R_d(x)$.

Observe that~\eqref{eq:tech-1} is simply a union bound over the contributions of all $d$ satisfying~\ref{tech-a}. By the principle of inclusion-exclusion, the sum of the first two terms in the inclusion-exclusion series are a lower bound on $R(x)$. The first term is simply the union bound that we have already computed. The second term is a sum over pairs $(d,d')$, both satisfying~\ref{tech-a}, where each term accounts for the contribution to $R(x)$ from those $n$ which are divisible by $\lcm(|d|,\!|d'|)$. Thus, abusing notation, we want to upper bound

\begin{equation}
    \label{eq:tech-2}
    \sum_{d,d' \t{ sat.\,\ref{tech-a}}} R_{d, d'}(x),
\end{equation}

where each term captures the contribution from $n$ dividing both $d$ and $d'$. As before, part of our proof involves showing that the sum
\[\sum_{d,d' \t{ sat.\,(a)}} \f{c_{dd'}}{\lcm(|d|,\!|d'|)}\]
converges for some $c_{dd'}$ which grows slowly as $\lcm(|d|,|d'|)$ goes to infinity. Here, we use~\cref{thm:lcm-sparse}, which tells us that the number of pairs of elements up to $x$ in a polynomially sparse set with least common multiple at most $x$ is $\ll x^{1-\kappa}$ for some $\kappa > 0$. By summation by parts, the sum in question converges. We use this to show that the second term in the inclusion-exclusion series is asymptotically negligible compared to the first term. Therefore, the constant we obtained from the union bound is actually the correct constant.

\subsection{Future work}
\label{subsec:future}

There are several natural extensions. The first concerns the generalization of our result to elliptic curves with good reduction everywhere and arbitrary $j$-invariant.

\begin{conjecture}
    \label{conj:irrational}
    Theorem~\ref{thm:rix-const} holds even after removing the constraint that $j$ is rational.
\end{conjecture}

To show this, we would want to show that the number of real and imaginary quadratic fields with discriminant of absolute value at most $x$ and over which there exists an elliptic curve with good reduction everywhere but no elliptic curve with good reduction everywhere and rational $j$-invariant is $o(\f{x}{\sqrt{\log x}})$. As mentioned, this conjecture is motivated by elliptic curve tables constructed by the first author~\cite{matschke} -- of the $955$ quadratic fields admitting an elliptic curve with good reduction everywhere, only $130$ of the fields do not admit such a curve that also has rational $j$-invariant. 
There are no known criteria as explicit as that of 
Setzer's result~\cite{setzer-crit} for identifying elliptic curves with good reduction everywhere and irrational $j$-invariant. However, the proof of 
Shavarevich's theorem~\cite{shafarevich62} (c.f.~Silverman~\cite[Thm.~IX.6.1]{silverman-book}) reduces the computation of elliptic curves over $K$ with good reduction everywhere (or more generally, outside a finite set of primes) to the computation of $\cO_K$-integral points on finitely many associated Mordell curves whose constant term lies in a Selmer-type group associated to $K$. This approach is used by the algorithm of~\cite{cremona-lingham}. In the case of elliptic curves with rational $j$-invariant, representatives of the relevant Selmer groups can be chosen to be rational integers, and in our language are the good $d$. Our proof in this paper thus has two steps: we show that the number of potential Selmer representatives is sparse (\cref{thm:red-ec} on how often $dy^2 = x^3-1728$ has a nontrivial integral point) and for each such potential representative we bound the number of quadratic fields $K$ for which it actually lies in the relevant Selmer-type group of $K$ (computing $R_d(x)$ and $I_d(x)$ for each good $d$). This compartmentalization appears to be more difficult in general as the former problem is no longer independent of the number field $K$ and because Selmer representatives can no longer be chosen to be rational integers. As such, it seems that a different approach may be needed in this setting.

Another angle of attack in the irrational $j$-invariant setting could come from the Shafarevich--Parshin construction~\cite{parshin72}, which reduces the computation of elliptic curves over $K$ with good reduction everywhere to the $S$-unit equation over $K$, where $S$ is the set of primes above~$2$. This approach was used in the aforementioned computation of~\cite{matschke}.

Another improvement would be removing the dependence on the $abc$-conjecture, which arises whenever we use~\cref{thm:abc-polyrad} to bound the range of $r$ for which it is possible for $\sqf(r^3-1728)$ to be $d$. As discussed in~\cref{sec:comp}, we can replace our dependence on the $abc$-conjecture with~\cref{assumption}, which implies both our result and Granville's~\cref{conj:granville-1}.

Lastly, it may be interesting to obtain criteria like that of Setzer for number fields of higher degrees. These criteria could then be used to derive statistical results for such families of number fields just as we have done in the quadratic case.

\subsection{Organization}

Most preliminary content, including proofs of elementary results and references to the literature, are in~\cref{sec:prelim}. We prove our bound on how often twists of some elliptic curves have an integral point, \cref{thm:red-ec}, in~\cref{sec:twists}. In~\cref{sec:asymp}, we apply Selberg-Delange theory and~\cref{thm:red-ec} to prove~\cref{thm:rix-asymp}, a version of~\cref{thm:rix-const} that is tight up to constants. In~\cref{sec:const}, we compute an upper bound on the leading constants in~\cref{thm:rix-const}. In~\cref{sec:lb}, we prove~\cref{thm:lcm-sparse} and use it to prove a matching lower bound on the leading constants of~\cref{thm:rix-const}, concluding the proof. In~\cref{sec:comp} we formulate an additional hypothesis based on~\cref{conj:granville-1} to obtain good numeric estimates for $c_R$ and $c_I$. 

\subsection{Acknowledgements}

The authors thank Ashwin Sah and Mehtaab Sawhney for an improvement to the upper bound in~\cref{thm:lcm-sparse}
and a matching lower bound.
In addition, the authors benefited from helpful conversations with 
Daniel Fiorilli, Andrew Granville, Michael Kural, and Melanie Matchett Wood. 
We thank the anonymous referee for their very useful comments.
The first author was supported by Boston University and by Simons Foundation grant \#550023.

\section{Preliminaries}
\label{sec:prelim}

\subsection{Notation}
\label{subsec:notation}

We write $\NN$ to denote the positive integers. In general, $p$ and $q$ will be used to denote primes and $\prod_p$, $\prod_q$, $\sum_p$, and $\sum_q$ denote products and sums over primes. For $p$ a prime, we write $|\argument|_p$ to denote the $p$-adic norm. Given a number field $K/\Q$, we write $\Delta_{K}$ to denote its absolute discriminant.

Let $n$ be an integer. We write $\omega(n)$ to denote the number of distinct prime factors of $n$. Generally, we will apply this in contexts where $n$ is squarefree, in which case $\omega(n)$ is simply the number of prime factors of $n$. 

If $t$ is the largest integer for which $t^2$ divides $n$ then we call $d := \f{n}{t^2}$ the \textdef{squarefree part} of $n$ and write $d = \sqf(n)$. The product of the distinct prime factors of $n$ is the \textdef{radical} of $n$, which we denote by $\rad(n)$. Note that the squarefree part of $n$ includes its sign but the radical does not. Given two positive integers $m$ and $n$, we write $(m,n)$ to denote their greatest common divisor.

Throughout, we use Vinogradov asymptotic notation.
If $f \ll g$ then $\limsup_{x \rightarrow \infty} \f{|f(x)|}{g(x)} < \infty$.
If $f \gg g$ then $\limsup_{x \rightarrow \infty} \f{|g(x)|}{f(x)} < \infty$.
If $f \ll g$ and $g \ll f$ then $f \asymp g$. 
We will also occasionally make use of Bachmann-Landau asymptotic notation to concisely describe error. The expression $b(x) = c(x) + o(g(x))$ means that  $\limsup_{x \rightarrow \infty} \f{b(x)-c(x)}{|g(x)|} = 0$. Similarly, we write $b(x) = c(x) + O(g(x))$ when $b(x) - c(x) \ll g(x)$. A subscript on any such notation -- for example, $\gg_{\eps}$ -- means that the implicit function or constant may depend on the subscript. We may sometimes combine both notations in expressions like $f(x) \ll (1+o_d(1))g(x)$ for some auxiliary variable $d$; this will be used if we wish to suppress the dependence on $d$ for brevity but will need to address it later in the paper. 
In this context, we write $o'_d(1)$ to denote an error term $\veps(d,x)$ such that 
\[\veps(d,x) \leq \paren{1-\f{\log |d|}{\log x}}^{-\nf{1}{2}}\paren{1+K'|d|^{1.002}\exp{\paren{-K\log^{\nf{1}{2}} \f{x}{|d|}}}+K'|d|^{0.001}\log^{-1}\f{x}{|d|}} - 1\] 
for some absolute constants $K$ and $K'$.
Note that for any large enough fixed $c$ (e.g.\ $c > 1$), $\eps(d,x)$ goes to zero in the regime 
where $x$ goes to infinity and $|d| \leq \log^c x$.

Given complex numbers $z$ and $a$, we define the power $z^a := e^{a\log z}$ with respect to the principal branch of the logarithm.

Throughout, we write $\leg{\argument}{\argument}$ to denote the Kronecker symbol.

We will often be concerned with Dirichlet characters of modulus $8$. To this end, it will help to fix a notation for characters of $(\Z/8\Z)^{\times} \cong (\Z/2\Z) \times (\Z/2\Z)$. Without loss of generality, let $3 \mod 8$ correspond to $(1,0) \in \Z/2\Z \times \Z/2\Z$, $5 \mod 8$ to $(0,1)$, and $7 \mod 8$ to $(1,1)$. We then define the characters $\boldsymbol{\chi}_{ij}$, for $i,j \in \{0,1\}$, to be nontrivial on the first (resp.\ second) component when $i$ (resp.\ $j$) is $1$. For example, $\boldsymbol{\chi}_{11}$ takes values $1,-1,-1,$ and $1$ at arguments which are $1,3,5,$ and $7 \mod 8$.  

\subsection{Some analytic facts}

\subsubsection{Summation by parts}
The following lemma will allow us to convert results on the sparseness of a subset of the natural numbers to results on the sums of reciprocals of elements of that subset.   
\begin{lemma}
    \label{lem:nat-to-harm}
    Let $\alpha \in (0,1)$ and $C$ be positive constants. Suppose that $f \colon \NN \rightarrow \NN$ is such that for all $x$,
        $\sum_{n = 1}^x f(n) \leq Cx^{1-\alpha}$.
    Then, if $\kappa +\alpha > 1$, 
    \begin{enumerate}[(i)]
        \item $\sum_{n=1}^{\infty} \f{f(n)}{n^{\kappa}}$ converges;\label{it:nat-to-harm-1}
        \item For $m \geq 2$, $\sum_{n = m}^{\infty} \f{f(n)}{n^\kappa} \leq \f{C\kappa m^{-(\kappa+\alpha-1)}}{\kappa+\alpha-1}$.\label{it:nat-to-harm-2}
    \end{enumerate}
\end{lemma}

\cref{lem:nat-to-harm} follows from a standard application of summation by parts or Stieljes integrals. 

\subsubsection{$L$-functions and convolutions}

Given two functions $a, b \colon \NN \rightarrow \C$ and associated Dirichlet series $F(s) = \sum_n a_nn^{-s}$ and $G(s) = \sum_n b_nn^{-s}$, their convolution is the formal Dirichlet series
\[F \otimes G(s) := \sum_n a_nb_nn^{-s}.\]

The convolution of $L$-series corresponds to the product of the coefficients in the same way that the product of $L$-series corresponds to the (Dirichlet) convolution of the coefficients.

Because we have only defined the convolution formally, we do need to consider the issue of convergence. However, in this paper we will only convolve $F$ and $G$ for which $|a_n| \leq 1$ and $|b_n| \leq 1$, and only ever require convergence of the convolution on $\Re(s) > 1$.

\subsubsection{Kronecker symbols}

The Kronecker symbol is a generalization of the Legendre symbol which allows composite inputs in the top and bottom entries and is - with some exceptions - multiplicative in both the top and the bottom entries. It is defined in most number theory texts (for example, page $39$ of~\cite{davenport}). Kronecker symbols are intimately related to real Dirichlet characters. The symbol $\leg{D}{\argument}$ is a real Dirichlet character when $D \not \equiv 3 \mod 4$, and every real Dirichlet character can be written as such a character. Furthermore, the primitive real Dirichlet characters are in $1-1$ correspondence with symbols $\leg{D}{\argument}$ when $D$ is a fundamental discriminant (i.e.\ when $D$ is the discriminant of a quadratic number field).

\subsection{The \texorpdfstring{$abc$}{abc}-conjecture and some consequences}

Recall the $abc$-conjecture.

\begin{conjecture}[Oesterl\'e~\cite{oesterle}, Masser~\cite{masser90}]
    \label{conj:abc}
    For every $\eps > 0$ there exists a constant $C_{\eps}$ such that for any given non-zero coprime integers $a,b,c$ with $a + b + c = 0$,
    \begin{equation*}
        \max(|a|,|b|,|c|) \leq C_{\eps}\rad(abc)^{1+\eps}.
    \end{equation*}
\end{conjecture}

\begin{theorem}[{Granville~\cite[Cor.~$1$]{granville-abc}}]
    \label{thm:abc-polyrad}
    Assume that the $abc$-conjecture (\cref{conj:abc}) is true. Suppose that $g(x) \in \Z[x]$ is separable. Then, for any $r \in \Z$, 
    \begin{equation*}
        \rad g(r) \gg_{\eps} |r|^{\deg g - 1 - \eps}.
    \end{equation*}
\end{theorem}

For a hyperelliptic curve over $\Q$ with integral model $C \colon y^2 = f(x)$, we write $C_d$ to denote its $d^{\t{th}}$ quadratic twist with model $dy^2 = f(x)$. We will make essential use of the following theorem.

\begin{theorem}[{Granville~\cite[Thm.~$1$(i)]{granville-twists}}]
    \label{thm:abc-rtd-bound}
    Assume that the $abc$-conjecture is true. Suppose that $f(x) \in \Z[x]$ is separable and let $C$ be the hyperelliptic curve with equation $y^2 = f(x)$. 
    If $\deg f \geq 3$ then the integral points $(r,t)$ on $C_d$ satisfy
    \begin{equation*}
        |r| \ll_{\eps} |d|^{\f{1}{\deg f - 2} + \eps}
    \end{equation*}
    and 
    \begin{equation*}
        |t| \ll_{\eps} |d|^{\f{1}{\deg f - 2} + \eps}
    \end{equation*}
    for every $\eps > 0$.
\end{theorem}

\begin{proof}
    Let $dt^2 = f(r)$. It follows from~\cref{thm:abc-polyrad} that, under the $abc$-conjecture
    \begin{equation*}
        |d|^{\f{1}{2}}|r|^{\f{\deg f}{2}} \gg_f |df(r)|^{\f{1}{2}} = |dt| \geq \rad(dt) = \rad(f(r)) \gg_{\eps} |r|^{\deg f - 1 - \eps}
    \end{equation*}
    and the first part of the result follows.
    For the second part, note that 
    \begin{equation*}
        |dt^2| = |f(r)| \ll_f |r|^{\deg f} \ll_{\eps} |d|^{\f{\deg f}{\deg f - 2} + \eps}.
        \qedhere
    \end{equation*}
\end{proof}

\subsection{Identifying quadratic fields with good reduction everywhere}
\label{subsec:setzer}

The following definition is the formal statement of~\ref{tech-a} from~\cref{subsec:tech}. 

\begin{definition}
We say that $d \in \Z$ is \textdef{good} if $d = \sqf(r^3-1728)$ for an element $r$ of the set
\begin{equation}
    \label{eq:good}
    \{r \in \Z \colon \t{if } r \equiv 0 \mod 2 \t{ then } r \equiv 0,4 \mod{16}; \t{ if } r \equiv 0\mod{3} \t{ then }n \equiv 12 \mod{27}\}.
\end{equation}
\end{definition}

Given $d \in \Z$ squarefree, we write 
\begin{equation}
    \eps_d :=
    \begin{cases}
        1 & d \equiv 1 \mod 4, \\
        -1 & \t{otherwise.} \\
    \end{cases}
\end{equation}

We may now state a criterion of Setzer~\cite{setzer-crit} which tells us when a quadratic field $\Q(\sqrt{m})$ admits an elliptic curve with good reduction everywhere and rational $j$-invariant. This will formalize~\ref{tech-b} and~\ref{tech-c} from~\cref{subsec:tech}. 

\begin{theorem}[Theorem 2.2 of~\cite{ctl}, correcting an error in Theorem 2 of~\cite{setzer-crit}]
    \label{thm:setzer-crit}
    Let $m$ be a squarefree integer. The field $\Q(\sqrt{m})$ admits an elliptic curve with good reduction everywhere and rational $j$-invariant if and only if the following conditions are satisfied for some integers $d$ and $n$ such that $d$ is good $($\eqref{eq:good}$)$ and $m = dn$.
    \begin{enumerate}[(i)]
        \item $\eps_dd$ is a quadratic residue modulo $n$; \label{it:setzer-1}
        \item $-\eps_d n$ is a quadratic residue modulo $d$; \label{it:setzer-2}
        \item If $d \equiv \pm 3 \mod 8$ then $m = dn \equiv 1 \mod 4$; \label{it:setzer-3}
        \item If $d$ is even then $n \equiv d+1 \mod 8$; \label{it:setzer-4}
        \item $m > 0$ if $\eps_dd < 0$. \label{it:setzer-5}
    \end{enumerate}
\end{theorem}

In the above, neither $d$ nor $n$ are restricted to being positive integers. Intuitively,~\ref{it:setzer-1} is~\ref{tech-b} from~\cref{subsec:tech}, and~\ref{it:setzer-2}-\ref{it:setzer-5} comprise~\ref{tech-c} from~\cref{subsec:tech}.

We will now define the primary quantities of interest. We write that an elliptic curve over a field $K/\Q$ has $\t{GRE}_{\Q}$ if it has good reduction at every prime of $\O_K$ and $j$-invariant in $\Q$.

\begin{definition}
    \label{def:ridx}
    Following~\cite{ctl}, we define for every positive $x$
    \begin{equation*}
        R(x) := \#\{m \colon 0 < \Delta_{\Q(\sqrt{m})} \leq x, \Q(\sqrt{m}) \t{ admits an elliptic curve with GRE}_{\Q}\}
    \end{equation*}
    and
    \begin{equation*}
        I(x) := \#\{m \colon -x \leq \Delta_{\Q(\sqrt{m})} < 0, \Q(\sqrt{m}) \t{ admits an elliptic curve with GRE}_{\Q}\}.
    \end{equation*}

    If $d$ is good, we also define
    \begin{equation*}
        R_d(x) := \#\{n \colon 0 < \Delta_{\Q(\sqrt{nd})} \leq x, \Q(\sqrt{nd}) \t{ admits an elliptic curve with GRE}_{\Q}\},
    \end{equation*}
    and
    \begin{equation*}
        I_d(x) := \#\{n \colon -x \leq \Delta_{\Q(\sqrt{nd})} < 0, \Q(\sqrt{nd}) \t{ admits an elliptic curve with GRE}_{\Q}\}.
    \end{equation*}
\end{definition}

Intuitively, $R_d(x)$ and $I_d(x)$ measure the contribution of $m$ which are divisible by $d$ to $R(x)$ and $I(x)$ respectively. Note that our definitions of $R_d(x)$ and $I_d(x)$ differ somewhat from those of Clemm and Trebat-Leder~\cite{ctl} as they write $R_d(x)$ (and $I_d(x)$ analogously) to count the number of $n$ up to $x$ such that $\Q(\sqrt{nd})$ admits an elliptic curve with $\t{GRE}_{\Q}$. 

\section{Twists of elliptic curves}
\label{sec:twists}

Let $E$ be an elliptic curve defined over $\Q$ with short Weierstrass equation $y^2 = f(x)$ for $f(x) \in \Z[x]$. We denote by $E_d$ the $d^{\t{th}}$ quadratic twist of $E$ with the integral model $dy^2 = f(x)$. Then, 
\begin{equation*}
    T_E(D) := \{d \in \Z \colon |d| \leq D, d \t{ is squarefree and } E_d \t{ has a nontrivial integral point}\}.
\end{equation*}

In this section, we will prove~\cref{thm:red-ec}, which we now recall. 

\redec

It will be useful to reformulate bounding $T_E(D)$ as bounding the number of integers with absolute value up to $D$ that arise as squarefree parts of $f(x)$. This is because the equation $dy^2 = f(x)$ has an integral point if and only if $d$ is the squarefree part of $f(r)$ for some integer $r$. The squarefree part function is multiplicative but not totally multiplicative. Nonetheless, the following lemma will let us write the squarefree part of a separable polynomial as the product of the squarefree parts of its factors if we ignore the valuations at finitely many primes.

\begin{lemma}
    \label{lem:squarefree-mult}
    Suppose that $f \in \Z[x]$ is separable and that it factors as $f = \prod_i g_i$ for $g_i \in \Z[x]$. Then, for all but finitely many $p$ and for all $r \in \Z$,
    \begin{equation*}
        |\sqf(f(r))|_p = \prod_i |\sqf(g_i(r))|_p
    \end{equation*}
\end{lemma}
\begin{proof}[Proof of~\cref{lem:squarefree-mult}]
    If two distinct factors $g_i$ and $g_j$ share a root modulo a prime then the resultant of $\bar{g}_i$ and $\bar{g}_j$ in $\Fp$ must be $0$, meaning that $p$ must divide the resultant of $g_i$ and $g_j$. By assumption, $g_i$ and $g_j$ have no shared roots over $\C$ and hence their resultant is a nonzero integer. This means that the reductions of $g_i$ and $g_j$ modulo a prime $p$ can only share a root for the finitely many primes $p$ that divide their resultant. Taking the union of these finite sets over all pairs of distinct $i$ and $j$, we see that for $p$ outside this finite union, if $p$ divides $g_i(r)$ for some $r$ then $p$ does not divide $g_j(r)$ for all $j \neq i$. This means that the factors $g_i(r)$ are ``coprime'' outside our finite set of primes, and the lemma then follows from the multiplicativity of the squarefree part function.
\end{proof}

\begin{proof}[Proof of \crefrestated{thm:red-ec}]
    There are four cases depending on how $f(x)$ factors. In the definition of each case, 
    a factorization of $f$ in $\Z[x]$ into irreducible factors is given.

    \begin{enumerate}
        \item $\boldsymbol{f(x) = (x+a_1)(x+a_2)(x+a_3)}$: This case follows directly from the proof of Theorem~$2$ in~\cite{granville-twists}.
        \item $\boldsymbol{f(x) = (x + a)(x^2 + b)}$: We have
            \begin{equation*}
              \f{1}{k(r)^2}\sqf(r+a)\sqf(r^2+b) = \sqf(f(r)),
            \end{equation*}
            where $k(r)$ is squarefree and divisible only by primes which divide both $r+a$ and $r^2+b$. By~\cref{lem:squarefree-mult},
            the set $S$ of primes $p$ for which there exists an $r$ such that $p$ divides both $r+a$ and $r^2+b$ is finite.
            Let $K := \prod_{p \in S} p$. Observe that $k(r) \leq K$. Then,
            \begin{equation*}
                \f{1}{K^2}\sqf(r+a)\sqf(r^2+b) \leq \sqf(f(r)),  
            \end{equation*}
            Choose an $\eps > 0$, suppose the constant from~\cref{thm:abc-rtd-bound} for this $\eps$ is $C_{\eps}$, and let
            \begin{equation}
                \label{eq:ohp}
                T'_E(D) := \{r \colon |r| \leq C_{\eps}D^{1+\eps}, |\sqf(r+a)\sqf(r^2+b)| \leq K^2D\}.
            \end{equation}
            Notice that $T_E(D)$ contains small values of $d$ while $T'_E(D)$ contains those $r$ for which $\sqf(f(r))$ is small.
            We claim that $|T_E(D)| \leq |T'_E(D)|$.
            This is because any $d \in T_E(D)$ is the squarefree part of $f(r)$ for some $r$. For each such $r$, 
            we have 
            \begin{equation}
              \label{eq:ohp2}
              |\sqf(r+a)\sqf(r^2+b)| = k(r)^2d \leq K^2D
            \end{equation}

            by~\eqref{eq:ohp} and hence $r \in T'_E(D)$.
            We only need to consider $|r| \leq C_{\eps}D^{1+\eps}$ by \cref{thm:abc-rtd-bound}.

            We will branch into two subcases. First, consider $r$ such that $|\sqf(r+a)| \leq D^{2\delta}$ for some $\delta > 0$ that we will select later. For each squarefree number $m \leq D^{2\delta}$, there are at most  $2\sqrt{\f{C_{\eps}D^{1+\eps}}{m}}$ values of $r$ such that $\sqf(r+a) = m$ and $|r| \leq C_{\eps}D^{1+\eps}$. Therefore, the contribution to $T'_E(D)$ from these $r$ is at most

            \begin{align*}
                2C_{\eps}^{\nf{1}{2}}\int_{1}^{D^{2\delta}} \paren{\f{D^{1+\eps}}{z}}^{\nf{1}{2}} dz & = 2C_{\eps}^{\nf{1}{2}}D^{\nf{1}{2}(1+\eps)}\int_{1}^{D^{2\delta}} z^{-\nf{1}{2}} dz \\
                & \leq C_1D^{\nf{1}{2}(1+\eps)+\delta}. \\
            \end{align*}
            for $C_1 \leq 4C_{\eps}^{\nf{1}{2}}$.

            Second,
            we will count $r$ contributing to $T'_E(D)$ for which $|\sqf(r+a)| > D^{2\delta}$.
            By~\eqref{eq:ohp2}, we require that $|\sqf(r^2+b)| \leq K^2D^{1- 2\delta}$.
            Fix an integer $m$ such that $|m| \leq K^2D^{1-2\delta}$.
            We will bound the number of $r$ for which $\sqf(r^2+b) = m$.
            For any such $r$,
            there is an $s \in \Z$ such that 
            \begin{equation*}
                r^2 - ms^2 = -b.
            \end{equation*}

            This is equivalent to requiring that $r+s\sqrt{m} \in \Z[\sqrt{m}]$ has norm $-b$.
            If $m$ is negative then $\Q(\sqrt{m})$ is imaginary quadratic and by Dirichlet's unit theorem there are only finitely many elements of $\Q(\sqrt{m})$ which have this norm.
            Hence, when $m$ is negative there are only finitely many values of $r$ for which $\sqf(r^2+b) = m$.

            If $m$ is positive, let $\eps_m$ be the unique fundamental unit that exceeds $1$ in the embedding of $\Q(\sqrt{m}) \hookrightarrow \R$ that sends $\sqrt{m} \mapsto \sqrt{m}$.
            Notice that after these choices $\eps_m > \f{1}{2}(1+\sqrt{m})$.
            This is because, by standard properties of the fundamental unit, $\eps_m = y+z\sqrt{m}$ for $y, z > 0$ and also $\eps_m$ is an element of $\cO_{\Q(\sqrt{m})}$.

            If $r+s\sqrt{m} \in \Z[\sqrt{m}]$ has norm $-b$ then any other element of $\Q(\sqrt{m})$ with this norm is, up to sign, a power of the fundamental unit $\eps_m$ times $r+s\sqrt{m}$. Suppose that $\alpha$ is an element (up to sign) of $\Z[\sqrt{m}]$ with norm $-b$ that has minimum archimedean absolute value. Note that there is a unique minimum archimedean absolute value because $\Z[\sqrt{m}]$ is discrete. Then, every element of $\Z[\sqrt{m}]$ with norm $-b$ is of form $\pm\alpha\eps_m^k$ for some $k \geq 0$. We know that $\eps_m > \f{1}{2}(1+\sqrt{m}) > 1.2$. We want to show that the $x$-coordinate of $\alpha \eps_m^k$ grows exponentially in $k$, as this will imply that we cannot have too many values of $r$ such that $\sqf(r^2+b) = m$. Writing $|\argument|$ for the archimedean absolute value, we see that if $r + s\sqrt{m} = \pm \alpha \eps_m^k$ then, 
            \begin{align*}
                |\alpha \eps_m^k| & = |r+s\sqrt{m}| \\
                & \leq |r|+|s\sqrt{m}| \\
                & = |r|+\sqrt{r^2+b} \\
                & \leq 2|r|+\sqrt{b}.
            \end{align*}
            Rearranging, this means that 
            \begin{equation*}
                \f{\alpha\eps_m^k-\sqrt{b}}{2} \leq |r|.
            \end{equation*}

            Because we require $|r| < C_{\eps}D^{1+\eps}$, we conclude that it is sufficient to look at 
            \begin{equation*}
                k < \f{\log\f{2C_{\eps}D^{1+\eps}+\sqrt{b}}{\alpha}}{\log \eps_m} \leq C_{2,\eps}\log D
            \end{equation*}
            for some positive $C_{2, \eps}$, using in the second inequality that $\eps_m$ is lower bounded by a constant independent of $m$.
            The point is that the number of possible $r$ that yield any given value of $\sqf(r^2+b)$ is small.
            Thus, the total number of possible $r$ in this subcase is at most $C_{2,\eps}K^2D^{1 - 2\delta}\log D$.

            Putting the two subcases together, our total count is $C_{1, \eps}D^{\f{1}{2}(1+\eps)+\delta} + K^2C_{2,\eps}D^{1-2\delta}\log D$.
            Taking $\delta = \f{1}{6}$, we can make our overall upper bound $(C_{1,\eps}+K^2C_{2,\eps})D^{\f{2}{3} + \eps}\ll_{\eps} D^{\f{2}{3} + \eps}$,
            since $K$ depends only on $f$.

        \item $\boldsymbol{f(x) = (x+a)(x^2+b_1x+b_2)}$: We will reduce to the previous case via a pair of coordinate transformations which map integral points to integral points and work for any twist of $y^2 = f(x)$. Start with the equation $dy^2 = f(x)$. Define $x' := 2^2x$, $y' := 2^3y$, $a' := 2^2a$, $b'_1 := 2^2b_1$, and $b'_2 := 2^4b_2$. We have
            \[dy'^2 = (x'+a')(x'^2+b'_1x'+b'_2),\]
            which looks the same as before but now $b'_1$ is even. Thus, we may complete the square, taking $x'' := x'+2^{-1}b'_1$, $y'' := y'$, $a'' := a'-2^{-1}b'_1$, and $b'' := b'_2 - 2^{-2}b'_1$ to obtain
            \[dy''^2 = (x''+a'')(x''^2+b'').\]
            Overall, we have $x'' = 2^2x + 2b_1$, $y'' = 2^3y$ and hence integral points map to integral points. Because the transformation is independent of $d$, we may apply the argument of the previous case to $f(x'') = (x''+a'')(x''^2+b'')$ to obtain an upper bound on the number of twists of our original equation with integral points.

        \item $\boldsymbol{f}$ \textbf{is not monic}: Suppose $f(x) = \sum_{i=0}^3 f_ix^i$.
          Consider $\tilde{x} := f_3x$, $\tilde{y} := f_3y$, and $\tilde{f}(\tilde{x}) := \tilde{x}^3 + f_2\tilde{x}^2 + f_1f_3\tilde{x} + f_0f_3^2$.
          Then, $\tilde{y}^2 = \tilde{f}(\tilde{x})$ and $\tilde{f}$ is monic.
          Furthermore, $\tilde{f}$ is reducible over $\Q$ because $f$ is.
          Thus, by Gauss' Lemma $\tilde{f}$ is reducible over $\Z$ and we may pass to one of the previous cases as appropriate because integral points are mapped to integral points. \qedhere
    \end{enumerate}
\end{proof}
\section{The Asymptotics of \texorpdfstring{$R(x)$}{R(x)} and \texorpdfstring{$I(x)$}{I(x)} up to Constants}
\label{sec:asymp}

Our main theorem in this section will be a version of~\cref{thm:rix-const} which is correct up to constants. 

\begin{theorem}
    \label{thm:rix-asymp}
    Assume that the $abc$-conjecture is true. Then, 
    \begin{equation*}
        R(x) \asymp \f{x}{\sqrt{\log x}}\quad  \t{and} \quad I(x) \asymp \f{x}{\sqrt{\log x}}.
    \end{equation*}
\end{theorem}

Because Clemm and Trebat-Leder proved the lower bound (\cref{thm:ctl}), it is sufficient for us to prove the upper bound. Throughout this section, $d$ will denote a squarefree integer. 
We will prove~\cref{thm:rix-asymp} by first proving the following lemma. 

Recall the definitions of $R_d(x)$ and $I_d(x)$ from \cref{def:ridx}, as well as the definition of $o'_d(1)$ from~\cref{subsec:notation}.
\begin{lemma}
    \label{lem:ridx-asymp}
    Let $d$ be good. Then, when $x \geq 3|d|$ we have
    \begin{equation*}
        R_d(x) \ll \f{(1+o'_d(1))c_dx}{|d|\sqrt{\log x}} \quad\t{and}\quad I_d(x) \ll \f{(1+o'_d(1))c_dx}{|d|\sqrt{\log x}},
    \end{equation*}
    where the implicit constant is absolute (independent of $d$) and $c_d$ is as defined in~\cref{cor:ridx-asymp}.
\end{lemma}

We will then pass from~\cref{lem:ridx-asymp} to~\cref{thm:rix-asymp} by summing $R_d(x)$ over good $d$ with small absolute value
and using a different bound for $d$ with large absolute value. By~\crefrestated{thm:red-ec},
the set of good $d$ is very sparse and therefore the asymptotic dependence on $x$ stays the same.
Notice that~\cref{lem:ridx-asymp} is unconditional -- we depend on the $abc$-conjecture in~\cref{thm:rix-asymp} only to prove~\crefrestated{thm:red-ec}.  
We keep track of the $o'_d(1)$ error
in~\cref{lem:ridx-asymp} in order to ensure that its dependence on $d$ is mild enough that summing over the set of good $d$ does not change the dependence on $x$ by more than a constant.

Throughout this section, we will compute all our bounds while pretending that our quadratic fields $\Q(\sqrt{m})$ 
(in both the real and imaginary settings) are ordered by $|m|$ rather than by $\Delta_{\Q(\sqrt{m})}$. 
This will not change any bound by more than a constant and thus is irrelevant to~\cref{lem:ridx-asymp} and~\cref{thm:rix-asymp}.

To prove~\cref{lem:ridx-asymp}, we will count those $n$ which satisfy~\cref{thm:setzer-crit}\ref{it:setzer-1}. This means that we want to count squarefree $n$ which are divisible only by primes $q$ for which $\eps_dd$ is a nonzero square modulo $q$. We do not actually need the squarefree condition to reach~\cref{thm:rix-asymp} but addressing it now will save us from repeating this work in~\cref{sec:const}. 
\begin{remark}[Difficulties of the sieve]
    Sieve theory provides a general method for counting the numbers up to $x$ which are divisible only by a subset $S$ of the primes.
    When $S$ has positive natural density $\alpha$, sieve theory gives us an upper bound of $O(x\log^{\alpha-1}x)$ for this number
    -- this follows from summation by parts and standard results.\footnote{See for example Theorem A.1 of~\cite{cribro}.} 
    However, in our setting $S$ is the set of primes satisfying~\cref{thm:setzer-crit}\ref{it:setzer-1}.
    By Chebotarev's density theorem applied to the multiquadratic extension $\Q(\sqrt{p_1}, \dots, \sqrt{p_r})$,
    one can check that $S$ has natural density $\f{1}{2}+ o_d(1)$ and our sieve-theoretic upper bound on $R_d(x)$ and $I_d(x)$ would end up being \[\ll \f{x}{d\log^{\nicefrac{1}{2}+ o_d(1)} x}.\] 
    This is insufficient to prove~\cref{lem:ridx-asymp}. 
    Furthermore, without the generalized Riemann hypothesis (or some other means of bounding Siegel zeroes) we have no control over the dependence on $d$ of the
    little-$o$ term that arises from Chebotarev's density theorem and thus cannot easily bound the sum of these answers over all good $d$.
    One can do better by exploiting the additional structure that the set $S$ has.
    Indeed, this is the idea behind the theorem of Serre~\cite{serre} (whose proof Serre attributed to Raikov, Wintner, and Delange)
    which Clemm and Trebat-Leder used in their lower bound.\footnote{Sieve-theoretic lower bounds which make much weaker assumptions on $S$ exist~\cite{gran-kou-mat, mat-shao} 
    but would have issues similar to those of the sieve-theoretic upper bounds were we to use them for our application.} 
    By using Selberg-Delange theory (\cref{thm:selberg-delange}) -- which generalizes some of the ideas used in~\cite{serre} -- we can obtain the correct dependence on $x$, 
    including the coefficient of $\f{x}{\sqrt{\log x}}$, and bypass the issue of computing the density of $S$ altogether. 
\end{remark}

\begin{definition}
    We define $\chi_d \colon \NN \rightarrow \C$ as the totally multiplicative function which on primes $q$ is
    \[\chi_d(q) := \begin{cases}
        \leg{\eps_dd}{q} & q \t{ odd}, \\
        \quad 1 & q=2 \t{ and } d \t{ odd}, \\
        \quad 0 & q=2 \t{ and } d \t{ even}, \\
    \end{cases}\]
    where $\leg{\argument}{\argument}$ denotes the Kronecker symbol. 
\end{definition}

Observe that unless $d \equiv \pm 3 \mod 8$ 
(and hence $\eps_dd \equiv 5 \mod 8$),
$\chi_d$ agrees with the Kronecker symbol $\leg{\eps_dd}{\argument}$ at all $n \in \NN$.
If $d \equiv \pm 3 \mod 8$, we define $\chi_d(2) = 1$ but $\leg{\eps_dd}{2} = -1$. By the definition of $\eps_d$,
we see that $\eps_dd \not \equiv 3 \mod 4$ and hence that $\leg{\eps_dd}{\argument}$ is necessarily a quadratic Dirichlet character. Furthermore, this character will always be primitive, as either $\eps_dd \equiv 1 \mod 4$, or $\eps_dd \equiv 2 \mod 4$ and $\leg{\eps_dd}{\argument} = \leg{4\eps_dd}{\argument}$. As such, $L(s,\chi_d)$, the $L$-function associated to $\chi_d$, is closely related to the $L$-function of some primitive quadratic Dirichlet character.

The following lemma captures some features of $\chi_d$ which we will use repeatedly throughout the paper. Henceforth, we write $m_d$ to denote the modulus of the Dirichlet character $\leg{\eps_dd}{\argument}$.
Recall from~\cref{subsec:notation} that $\boldsymbol{\chi}$ is used to refer to the quadratic characters modulo $8$. In particular, $\boldsymbol{\chi}_{11}$ should not be confused with $\chi_d$ for $d = 11$, and the latter will never appear in this paper. 
\begin{lemma}
    \label{lem:chid}
    \leavevmode
    \begin{enumerate}[(i)]
        \item Let $d$ be odd, and write $d = \pm p_1\cdots p_r$. Then, for $n$ odd, we have 
        \[\chi_d(n) = \prod_{1 \leq i \leq r} \leg{n}{p_i} \quad\t{and}\quad \chi_d(2n) = \chi_d(n).\]
        Let $d$ be even, and write $d = \pm 2p_1p_2\cdots p_r$. Then, we have 
        \[\chi_d(n) = \begin{cases}
            \boldsymbol{\chi}_{01}(n) \prod_{1 \leq i \leq r} \leg{n}{p_i} & d \equiv 2 \mod 8, \\
            \boldsymbol{\chi}_{11}(n) \prod_{1 \leq i \leq r} \leg{n}{p_i} & d \equiv 6 \mod 8. \\
        \end{cases}\] \label{it:chid-1}
        \item 
        We have
        \[L(s, \chi_d) = \begin{cases}
        L(s,\leg{\eps_dd}{\argument})\f{1+2^{-s}}{1-2^{-s}} & d \equiv \pm 3 \mod 8, \\
        L(s, \leg{\eps_dd}{\argument}) & \t{otherwise.} \\
        \end{cases} \] \label{it:chid-2}
        \item We have
        \[m_d = \begin{cases}
            |d| & d \t{ odd}, \\
            4|d| & d \t{ even}. 
        \end{cases}\] \label{it:chid-3}
    \end{enumerate}
\end{lemma}

\begin{proof}[Proof sketch]
\leavevmode
\begin{enumerate}[(i)]
    \item This is a standard application of quadratic reciprocity. 
    \item If $d \not \equiv \pm 3 \mod 8$, $\chi_d = \leg{\eps_dd}{\argument}$. Otherwise, $\chi_d$ is multiplicative and agrees with $\leg{\eps_dd}{\argument}$ except at $2$, where it has value $1$ instead of $-1$. 
    \item This is a standard property of the Kronecker symbol. \qedhere
\end{enumerate}
\end{proof}

\begin{definition}
    \label{def:ns}
    We define $N_{d}(y)$ to be the set of squarefree integers between $1$ and $y$ which are divisible only by primes $q$ for which $\chi_d(q) = 1$. 
\end{definition}

We will ultimately use $N_d(\f{x}{|d|})$ to upper bound $R_d(x)$ and $I_d(x)$. Let
\begin{equation*}
    a_n := \begin{cases} 
        1 & n \t{ is squarefree and divisible only by primes $q$ s.t. } \chi_d(q) = 1 \\
        0 & \t{otherwise.} \\
    \end{cases}
\end{equation*}

We will obtain~\cref{lem:ridx-asymp} by deriving the asymptotic behavior of $N_d(y) = \sum_{n \leq y} a_n$. We can access the latter by studying the Dirichlet series defined on $\Re(s) > 1$ by
\begin{align*}
    F_d(s) := \sum_{n \geq 1} a_nn^{-s}.
\end{align*}

This next lemma is more general than we need in this section but we will use its full generality later. 

\begin{lemma}
    \label{lem:ridx-factor-conv}
    If $L(s,\chi)$ is a Dirichlet $L$-series, we have
    \[F_d(s) \otimes L(s,\chi) = C(s)L(s,\chi)^{\nf{1}{2}}L(s,\chi_d\chi)^{\nf{1}{2}},\]
    where
    \[C(s) := L(2s,\chi^2)^{-\nf{1}{2}}\prod_{q \in S} (1+\chi(q)q^{-s})^{-\nf{1}{2}}\prod_{q:\chi_d(q)=1} (1-\chi^2(q)q^{-2s})^{\nf{1}{2}}\] 
    extends to a holomorphic function on $\Re(s) > \f{1}{2}$. 
\end{lemma}

\begin{proof}
    By the definition of convolution, 
    \[F_d(s) \otimes L(s,\chi) = \sum_{n \geq 1} a_n\chi(n)n^{-s}.\]
    This is holomorphic on $\Re(s) > 1$. Because $a_n$ is $0$ whenever $n$ is not squarefree, we have the Euler product
    \begin{align*}
        {F_d}(s) \otimes L(s,\chi)  & = \prod_q (1+\chi(q)a_qq^{-s}) \\
        & = \prod_q (1-\chi(q)a_qq^{-s})^{-1}\prod_q (1-\chi^2(q)a^2_qq^{-2s}). \\
      \end{align*}
    Because $a_q \in \{0, 1\}$ for all $q$, we have
    \begin{align*}
        {F_d}(s) \otimes L(s,\chi)  & = \prod_q (1-\chi(q)q^{-s})^{-\nf{1}{2}}(1-\chi(q)\chi_d(q)q^{-s})^{-\nf{1}{2}} \\
        & \hspace{1cm}\times \prod_{q:\chi_d(q)=\t{-}1} (1-\chi^2(q)q^{-2s})^{\nf{1}{2}}\prod_{q|m_d} (1-\chi(q)q^{-s})^{\nf{1}{2}}\prod_q (1-\chi^2(q)a^2_qq^{-2s}). \\
    \end{align*}
    The first two terms are the local factors $L(s,\chi)^{\nf{1}{2}}$ and $L(s,\chi\chi_d)^{\nf{1}{2}}$. For the other three, we have
    \begin{align*}
        & \prod_{q:\chi_d(q)=\t{-}1} (1-\chi^2(q)q^{-2s})^{\nf{1}{2}}\prod_{q|m_d} (1-\chi(q)q^{-s})^{\nf{1}{2}}\prod_q (1-\chi^2(q)a^2_qq^{-2s}) \\
        = \ & \prod_{q:\chi_d(q)=\t{-}1} (1-\chi^2(q)q^{-2s})^{\nf{1}{2}}\prod_{q\mid m_d}(1-\chi(q)q^{-s})^{\nf{1}{2}}\prod_{q:\chi_d(q)=1} (1-\chi^2(q)q^{-2s}) \\
        = \ & \prod_q (1-\chi^2(q)q^{-2s})^{\nf{1}{2}} \prod_{q \mid m_d} (1+\chi(q)q^{-s})^{-\nf{1}{2}}\prod_{q:\chi_d(q)=1} (1-\chi^2(q)q^{-2s})^{\nf{1}{2}} \\
        = \ & L(2s,\chi^2)^{-\nf{1}{2}}\prod_{q \mid m_d} (1+\chi(q)q^{-s})^{-\nf{1}{2}}\prod_{q:\chi_d(q)=1} (1-\chi^2(q)q^{-2s})^{\nf{1}{2}},
    \end{align*}
    which is the expression for $C(s)$ in the lemma. This product converges absolutely and locally uniformly on $\Re(s) > \f{1}{2}$ so we have holomorphicity on this region as desired. 
\end{proof}

In this section, we will use the following simple corollary of~\cref{lem:ridx-factor-conv}. 
\begin{corollary}
    \label{cor:ridx-factor-sqf}
    Let $F_d(s)$ be as defined in~\cref{lem:ridx-factor-conv}. Then,
    \[F_d(s) = C_d(s)\zeta(s)^{\nicefrac{1}{2}}L(s, \chi_d)^{\nicefrac{1}{2}},\]
    where 
    \[C_d(s) := \zeta(2s)^{-\nicefrac{1}{2}}\prod_{q \mid m_d} (1+q^{-s})^{-\nicefrac{1}{2}}\prod_{q:\chi_d(q)=1} (1-q^{-2s})^{\nicefrac{1}{2}}\]
    extends to a holomorphic function on $\Re(s) > \f{1}{2}$. 
\end{corollary}

\begin{proof}
    Apply~\cref{lem:ridx-factor-conv} with $\chi$ as the trivial character.  
\end{proof}

To pass from the factorization of a Dirichlet series to a bound on the sum of its coefficients, we will use a Tauberian theorem given by the main theorem of Selberg-Delange theory. 

\begin{theorem}[Special case of Chapter II.5.3 Theorem 3 from \cite{tenenbaum}]
    \label{thm:selberg-delange}
    Suppose that $F(s) := \sum_{n} h_nn^{-s}$ is a Dirichlet series such that $h_n \geq 0$ for all $n$. Suppose further that for some $\rho \in \C$, 
    the function
    ${G(s) := F(s)\zeta(s)^{-\rho}}$ can be analytically continued to a holomorphic function on the region \linebreak $\Re(s) \geq 1-\beta$ for some $\beta > 0$, and in this region satisfies the bound 
    \begin{equation}
        \label{eq:selberg-delange-1}
        |G(s)| \leq M(1+|\Im(s)|)^{1-\delta}
    \end{equation}
    for some $M > 0$ and $0 < \delta \leq 1$. Then, if $x \geq 3$ we have
    \[\sum_{n \leq x} h_n = \f{x}{\log^{1-\rho} x}\paren{\f{G(1)}{\Gamma(\rho)} + O\big(Me^{-K\sqrt{\log x}} + \log^{-1} x\big)}\]
    where $K$ and the implicit constant are absolute.
\end{theorem}

With~\cref{thm:selberg-delange} in hand we may now pass from properties of $F_d(s)$ to the asymptotics of $N_d(y)$.

\begin{corollary}
    \label{cor:ridx-asymp}
    We have
    \begin{equation}
        \label{eq:ridx-asymp-1}
        |N_d(y)| = \sum_{n \leq y} a_n = \f{(c_d+\veps_1(d,y))y}{\sqrt{\log y}},
    \end{equation}
    where 
    \begin{equation}
        \label{eq:ridx-asymp-2}
        c_d = \f{1}{\Gamma(\nicefrac{1}{2})\zeta(2)^{\nicefrac{1}{2}}}L(1,\chi_d)^{\nicefrac{1}{2}}\prod_{q \mid m_d} \paren{1+q^{-1}}^{-\nicefrac{1}{2}}\prod_{q:\chi_d(q) = 1} \paren{1-q^{-2}}^{\nicefrac{1}{2}}.
    \end{equation}
    and 
    \[\veps_1(d,y) \ll |d|^{1.001}e^{-K\sqrt{\log y}}+\log^{-1} y.\]

    If $H(s) = A(s)F_d(s) = \sum_n h_nn^{-s}$ on $\Re(s) > 1$ for some $A(s)$ which extends to a bounded holomorphic function on $\Re(s) > \f{1}{2}$,
    \[\absBigg{\sum_{n \leq y} h_n} = \f{(A(1)c_d+\veps_2(d,y))y}{\sqrt{\log y}}.\] 
    where again
    \[\veps_2(d,y) \ll |d|^{1.001}e^{-K\sqrt{\log y}}+\log^{-1} y.\]
\end{corollary}
 
We will need a simple bound in the proof of~\cref{cor:ridx-asymp}.

\begin{proposition}[Folklore, e.g.\ \cite{granville-soundararajan}]
    \label{prop:l-bound}
    Let $\chi$ be a nonprincipal Dirichlet character of modulus $m$. 
    Then, $m^{-\eps} \ll_{\eps} |L(1,\chi)| \ll \log m$ for $\eps > 0$. 
\end{proposition}

\begin{proposition}[Phragm\'en-Lindel\"of in a strip~\cite{iwaniec-kowalski}]
  \label{prop:phragmen-lindelof}
  Let $\Omega := \{s \in \C \colon a \leq \Re(s) \leq b\}$.
  For $t \in \R$, suppose that
  \begin{itemize}
    \item $f$ is holomorphic on an open neighborhood of $\Omega$,
    \item $|f(s)| \ll e^{|s|}$ on $\Omega$,
    \item $f(a+it) \leq M_a(1+|t|)^{\alpha}$, and
    \item $f(b+it) \leq M_b(1+|t|)^{\beta}$.
  \end{itemize}

  Let $\ell(x) := \f{b-x}{b-a}$. Then, 
\[|f(s)| \leq M_a^{\ell(\Re(s))}M_b^{1-\ell(\Re(s))}(1+|\Im(s)|)^{\alpha\ell(\Re(s))+\beta(1-\ell(\Re(s)))}\]
  for all $s \in \Omega$. 
\end{proposition}

\begin{lemma}
  \label{lem:pl-l-func}
  Let $\chi$ be a primitive Dirichlet character of modulus $m$.
  On the region $\Re(s) \geq \f{1}{2}$, we have
  \[|L(s,\chi)| \ll m^{\nf{1}{4}}(1+|\Im(s)|)^{\nf{1}{4}},\]
  where the implicit constant is absolute (independent of $m$ and $\chi$).
\end{lemma}

\begin{proof}
  Let $\Omega := \{s \in \C \colon \Re(s) > \f{1}{2}\}$.
  The convexity bound (see e.g. Theorem 5.23 of~\cite{iwaniec-kowalski}) 
  for Dirichlet $L$-functions
  tells us that when $\Re(s) = \f{1}{2}$ (i.e.\ on $\bar{\Omega}$)
  \begin{equation}
    \label{eq:convexity}
    |L(s,\chi)| \leq Km^{\nf{1}{4}}|s|^{\nf{1}{4}}
  \end{equation}
  for some absolute constant $K > 0$. 
  When $\Re(s) > 1$, we know that
  \[|L(s,\chi)| \leq \zeta(\Re(s)) \leq \f{\Re(s)}{\Re(s)-1}.\]
  For any $s \in \Omega$, we can thus apply~\cref{prop:phragmen-lindelof} 
  to the strip with $a = \f{1}{2}$ and $b = \Re(s)+1$, 
  $\alpha = \f{1}{4}$, $\beta = 0$, $M_a = Km^{\nf{1}{4}}$ (from~\eqref{eq:convexity}),
  and $M_b = \f{\Re(s)+1}{\Re(s)}$.~\cref{prop:phragmen-lindelof} then tells 
  us that 
  \begin{align*}
    |L(s,\chi)| & \leq (Km^{\nf{1}{4}})^{\ell(\Re(s))}\paren{\f{\Re(s)+1}{\Re(s)}}^{1-\ell(\Re(s))}(1+|\Im(s)|)^{\nf{1}{4}\ell(\Re(s))} \\
                & \leq 3\max(K,1)m^{\nf{1}{4}}(1+|\Im(s)|)^{\nf{1}{4}} \\
  \end{align*}
  for every $s \in \Omega$, as desired.
\end{proof}

\begin{theorem}[Robin~\cite{robin}]
    \label{thm:robin}
    For $n \in \NN$, let $\sigma(n)$ denote the sum of the positive divisors of $n$. Then, for all $n \geq 3$,
    \[\sigma(n) < e^{\gamma}n\log\log n+\f{0.6483n}{\log\log n},\]
    where $\gamma$ is the Euler-Mascheroni constant. 
\end{theorem}

\begin{corollary}
    \label{cor:robin}
    For $n \in \NN$ squarefree and $s \in \C$ with $\Re(s) \geq 1-\beta$ for $\beta \in [0,1)$,  
    \[\absBigg{\prod_{q | n} (1-q^{-s})} \leq n^{\beta}\paren{e^{\gamma}\log\log n+\f{0.6483}{\log\log n}}.\]
\end{corollary}
\begin{proof}
        \[\absBigg{\prod_{q | n} (1-q^{-s})} \leq \prod_{q | n} (1+q^{-(1-\beta)}) 
         = \prod_{q | n} \f{q^{1-\beta} + 1}{q^{1-\beta}}
         \leq \prod_{q | n} \f{q + 1}{q^{1-\beta}}
         \leq \f{\sigma(n)}{n^{1-\beta}},\]
    and applying~\cref{thm:robin} yields the result. 
\end{proof}

\begin{proof}[Proof of~\cref{cor:ridx-asymp}]
    Recall the factorization of $F_d(s)$ that we obtained from~\cref{cor:ridx-factor-sqf}.
    We want to apply~\cref{thm:selberg-delange} with $\rho = \frac{1}{2}, \delta < \f{1}{2}$ and 
    \begin{equation}
        G(s) = L(s,\chi_d)^{\nicefrac{1}{2}}C_d(s),
    \end{equation}
    since we know that $L(s,\chi_d)$ can be analytically continued to $\C$ and $C_d$ is holomorphic on the region $\Re(s) > \frac{1}{2}$.
    Fix $\beta < 0.0005$.
    We need to bound $G(s)$ on the region $\Re(s) > 1-\beta$.
    By~\cref{lem:chid} and~\cref{lem:pl-l-func} we have
    \[|L(s,\chi_d)| \ll |L(s,\leg{\eps_dd}{\argument})| \ll m_d^{\nf{1}{4}}(1+|\Im(s)|)^{\nf{1}{4}}.\]
    By~\cref{cor:robin},
    \[\absBigg{\prod_{q | m_d} (1-q^{-s})^{\nf{1}{2}}} \ll m_d^{\nf{\beta}{2}}\paren{\log\log m_d + \f{1}{\log\log m_d}}^{\nf{1}{2}}.\]

    Lastly, observe that \[\abs{\prod_{q:\chi_d(q) = \t{-}1} (1-q^{-2s})^{\nicefrac{1}{2}}} \leq |\zeta(2(1-\beta))^{\nicefrac{1}{2}}|.\]
 
    This means that~\cref{thm:selberg-delange} can be applied with $M \ll |d|^{1.001+\nf{\beta}{2}}$, giving us~\eqref{eq:ridx-asymp-1} with constant 
    \[c_d = \frac{G(1)}{\Gamma(\nicefrac{1}{2})}\]
    and error
    \[\veps_1(d,y) \ll |d|^{1.001}e^{-K\sqrt{\log y}}+\log^{-1} y. \qedhere\]
\end{proof}

\begin{lemma}
  \label{lem:cd-lb}
  \[c_d \gg |d|^{-0.001}.\]
\end{lemma}

\begin{proof}
  Recall from~\cref{cor:ridx-factor-sqf} that  
  \[c_d := \f{L(1,\chi_d)^{\nicefrac{1}{2}}}{\Gamma(\nicefrac{1}{2})\zeta(2)^{\nicefrac{1}{2}}} \
  \prod_{q \mid m_d} \paren{1+q^{-1}}^{-\nicefrac{1}{2}}\prod_{q:\chi_d(q) = 1} \ 
  \paren{1-q^{-2}}^{\nicefrac{1}{2}}.\]

  We have the lower bounds
  \begin{equation}
    \label{eq:oof-1}
    \prod_{q:\chi_d(q) = 1} \ 
    \paren{1-q^{-2}}^{\nicefrac{1}{2}} \geq \zeta(2)^{-\nicefrac{1}{2}}
  \end{equation}

  and 
  \begin{equation}
    \label{eq:oof-2}
    \prod_{q \mid m_d} \paren{1+q^{-1}}^{-\nicefrac{1}{2}} \
    \geq e^{-\f{1}{2}\sum_{q \mid 4d} \f{1}{q}} \
    \geq e^{-\f{1}{2}(K + \log \log (\ceil{\log_2 |4d|}+1))} \
    \geq \paren{e^K\log (\log_2 |4d| + 2)}^{-\nf{1}{2}},
  \end{equation}
  where in the first inequality we have used that $4d$ has at
  most $\ceil{\log_2 |4d|}$ prime factors and and in the second inequality
  $K$ is some absolute constant and we have used that the sum of the reciprocals
  of the first $n$ primes is at most $\log \log (n+1) + 1$. 

  Applying~\cref{prop:l-bound}, \cref{eq:oof-1}, and \cref{eq:oof-2},
  for any $\eps' > \eps > 0$
  \[c_d \gg_\eps \f{1}{|d|^{\eps}(\log\log |d|)^{\nf{1}{2}}} \gg_{\eps'} |d|^{-\eps'}\]
  as desired. 
\end{proof}

\begin{proof}[Proof of~\cref{lem:ridx-asymp}]
    Recall that $R_d(x)$ and $I_d(x)$ are the number of $n$ in the appropriate intervals
    which satisfy all of the criteria in~\cref{thm:setzer-crit}.
    Because the lower bound is already known, all we need is an upper bound.
    We give the proof for $R_d(x)$ but because this criterion is independent of the sign of $n$ the same approach will work for $I_d(x)$. 
    It suffices to count those $n$ which satisfy just~\ref{thm:setzer-crit}\ref{it:setzer-1}.
    Observe that $R_d(x) \leq N_d(\f{x}{|d|})$ because we require that $|dn| \leq x$.

    Plugging this into~\cref{cor:ridx-asymp}, we have 
    \[R_d(x) \leq \f{(c_{d}+\veps_1(d,\f{x}{|d|}))x}{|d|\sqrt{\log x - \log |d|}}.\]
    We can replace $\sqrt{\log x - \log |d|}$ in the denominator with $\sqrt{\log x}$ by multiplying the 
    numerator and the denominator by
    \[\sqrt{\f{\log x}{\log x - \log |d|}} = \paren{1-\f{\log |d|}{\log x}}^{-\nf{1}{2}}.\]

    By~\cref{lem:cd-lb}, 
    \[R_d(x) \leq \f{\paren{c_{d}+\veps_1(d,\f{x}{|d|})}\paren{1-\f{\log |d|}{\log x}}^{-\nf{1}{2}}x}{|d|\sqrt{\log x}} \ll \f{(1+o'_d(1))c_dx}{|d|\sqrt{\log x}}\]
    as desired. 
\end{proof}

\begin{proof}[Proof of \cref{thm:rix-asymp}]
    We will give the proof for $R(x)$; the analogous proof works for $I(x)$. If $\Q(\sqrt{m})$ admits a curve with good reduction everywhere then some good $d$ must divide $m$. There could be multiple good $d$ that divide a given $m$, but at least as an upper bound we have
    \begin{equation*}
        R(x) \leq \sum_{\substack{d \t{ good} \\ |d| \leq x}} R_d(x)
    \end{equation*}

    Let $z := \log^{\nicefrac{3}{2}+\delta}(x)$ for $\delta > 0$.
    Per~\cref{lem:ridx-asymp}, we have that
    \begin{equation}
        \label{eq:rix-ub-1}
        \sum_{\substack{d \t{ good} \\ |d| \leq x}} R_d(x) \ll x\paren{\sum_{\substack{d \t{ good} \\ |d| \leq z}} \f{(1+o'_d(1))c_d}{|d|\sqrt{\log x}} + \sum_{\substack{d \t{ good} \\ |d| \geq z}} \f{1}{|d|}},
    \end{equation}
    where for the second sum we have used that we always have a bound of $\f{x}{|d|}$ on the number of natural numbers up to $x$ which are multiples of $|d|$,
    and in the first sum we have used that for large enough $x$ (say, $x > 1000$), $\f{x}{|d|} \geq \f{x}{z} \geq 3z \geq 3|d|$ so we are in the regime where
    ~\cref{lem:ridx-asymp} holds.

    \cref{thm:red-ec} tells us that the set of good $d$ is $\f{1}{3}$-polynomially sparse. Therefore, by~\cref{lem:nat-to-harm}\ref{it:nat-to-harm-2},

    \[\sum_{\substack{d \t{ good} \\ |d| \geq z}} \f{1}{|d|} \ll z^{-\nicefrac{1}{3}}.\]

    By our choice of $z$, the second term 
    on the right-hand-side of~\eqref{eq:rix-ub-1} will be negligible compared to the first.
    By the definition of $o'_d(1)$ in~\cref{subsec:notation}, the error term in the numerator of the first term on the right-hand-side of~\eqref{eq:rix-ub-1} goes to zero as $x$ goes to infinity when $|d| \leq z$.

    Putting everything together, we can rewrite~\eqref{eq:rix-ub-1} as 
    \begin{equation*}
        \sum_{\substack{d \t{ good} \\ |d| \leq x}} R_d(x) \ll \f{x}{\sqrt{\log x}}\sum_{\substack{d \t{ good} \\ |d| \leq z}} \f{c_d}{|d|}. 
    \end{equation*}

    Now we upper bound $c_d$. By applying~\cref{prop:l-bound} and~\cref{cor:robin} (with $\beta = 0$)  to~\eqref{eq:ridx-asymp-2}, we see that
    \[c_d \ll |d|^{0.001},\]
    meaning we need to bound the sum of $\f{1}{|d|^{.999}}$. Because $.999 + \f{1}{3} > 1$, applying~\cref{lem:nat-to-harm}\ref{it:nat-to-harm-1} with $\alpha = \f{1}{3}$ and $\kappa = .999$ tells us that the sum in question converges. Therefore,
    \begin{equation*}
        R(x) \ll \f{x}{\sqrt{\log x}}.
    \end{equation*}
    as desired.
\end{proof}

\section{An Upper Bound on the Constant}
\label{sec:const}
We are now ready to begin the proof of our main theorem, which we recall for convenience. 

\rixconst

In this section, we show that the expressions in~\crefrestated{thm:rix-const} are upper bounds on the correct values of $c_R$ and $c_I$.

The main lemma of this section gives a sharp version of~\cref{lem:ridx-asymp}. Once again recall the definition of $o'_d(1)$ from~\cref{subsec:notation}. 

\begin{lemma}
    \label{lem:ridx-const}
    Let $d$ be good. Then, 
    \begin{equation*}
        R_d(x) = \f{(1+o'_d(1))c_dc'_{d,R}x}{|d|2^{\omega(d)}\sqrt{\log x}}
        \quad\t{and}\quad
        I_d(x) = \f{(1+o'_d(1))c_dc'_{d,I}x}{|d|2^{\omega(d)}\sqrt{\log x}},
    \end{equation*}
    where $c_d$ is as defined in~\cref{cor:ridx-asymp},
    \[c'_{d,R} := \begin{cases}
        1 & d \equiv \pm 1 \mod 8 \\
        \f{2}{3} & d \equiv \pm 3 \mod 8 \\
        \f{1}{4} & d \equiv 2 \mod 8 \quad\t{and}\quad d > 0 \\
        \f{1}{4} & d \equiv 6 \mod 8 \\
        0 & \t{otherwise,}
    \end{cases}
    \]
    \[
    c'_{d,I} := \begin{cases}
        1 & d \equiv 1 \mod 8 \quad\t{and}\quad d > 0 \\
        1 & d \equiv 7 \mod 8 \quad\t{and}\quad d < 0 \\
        \f{2}{3} & d \equiv 3 \mod 8 \quad\t{and}\quad d < 0 \\
        \f{2}{3} & d \equiv 5 \mod 8 \quad\t{and}\quad d > 0 \\
        \f{1}{4} & d \equiv \pm 2 \mod 8 \quad\t{and}\quad d < 0 \\
        0 & \t{otherwise.}
    \end{cases}
    \]
\end{lemma}

\begin{remark}
    Notice that $c'_{d,R} \geq c'_{d,I}$ except when $d$ is a negative number congruent to $2 \mod 8$. 
    Per~\crefrestated{thm:rix-const}, 
    \[c_R - c_I = \sum_{d \t{ good}} \f{c_d(c'_{d,R}-c'_{d,I})}{|d|2^{\omega(d)}}.\]
    This provides strong evidence that $c_R > c_I$ under the $abc$-conjecture, but does not provide a proof without better control over the distribution of good $d$ and its correlation with $\f{c_d}{2^{\omega(d)}}$ than we are able to show. We prove $c_R > c_I$ in a different manner, assuming a different hypothesis, in~\cref{sec:comp}. 
\end{remark}

In the proof of~\cref{lem:ridx-asymp} it was sufficient to consider~\cref{thm:setzer-crit}\ref{it:setzer-1}. Now, we will need to consider all five conditions. Observe that while~\ref{thm:setzer-crit}\ref{it:setzer-1} imposes a condition on the primes that are allowed to divide~$n$,~\ref{it:setzer-2}-\ref{it:setzer-4} constrain the value of $n$ modulo the primes dividing $d$ and modulo $4$ or $8$. These thus correspond to~\ref{tech-c} in~\cref{subsec:tech}. As discussed then, the Dirichlet series (of the indicator function) for this property is not multiplicative, complicating any application of~\cref{thm:selberg-delange}. We address this by expressing the series as a linear combination of Dirichlet-$L$-series. Intuitively, a Tauberian theorem can be thought of as telling us the “rate of divergence” at a pole, and hence only terms which possess a pole at $s = 1$ will contribute to the overall asymptotic.

\begin{proof}
    We start with the bound for $R_d(x)$. If $d$ is odd we write $d = \pm p_1p_2 \dots p_r$ and if $d$ is even we write $d = 2d' = \pm 2p_1p_2\dots p_r$. We will count positive numbers $n$ such that $\Q(\sqrt{dn\sgn(d)})$ has discriminant with absolute value at most $x$ and admits an elliptic curve with good reduction everywhere and rational $j$-invariant. \textit{Note that this $n$ is different from the $n$ in~\cref{thm:setzer-crit} as it is always positive.}

    For any valid pair of $(d,n)$,~\cref{thm:setzer-crit} tells us that the following conditions on $n$ and $d$ hold:
    \begin{enumerate}[(a)]
        \item $n$ is coprime to $d$;\footnote{This condition is redundant in light of~\ref{it:const-1c} but we include it for clarity.}\label{it:const-1a}
        \item $n$ is squarefree;\label{it:const-1b}
        \item $\chi_d(q) = 1$ for every $q \mid n$;\label{it:const-1c}
        \item $\leg{n}{p} = \leg{-\sgn(d)\eps_d}{p}$ for every odd $p \mid d$.\label{it:const-1d}
    \end{enumerate}

    Let us check what we need for these conditions to be compatible with one another. Assume that~\ref{it:const-1a} and~\ref{it:const-1b} hold for some $d$ and $n$. By~\cref{lem:chid},~\ref{it:const-1c} tells us the value of the product
    \[\prod_{p \mid d} \leg{q}{p}\]
    for each $q | n$.
    Similarly,~\ref{it:const-1d} tells us something about 
    \begin{equation*}
         \prod_{q | n} \leg{q}{p} = \leg{n}{p}
    \end{equation*}
    for each $p | d$. Let $M$ be a matrix with rows and columns indexed by $p_i | d$ and $q_j | n$ respectively and with entries $M_{ij} := \leg{q_j}{p_i}$. Then,~\ref{it:const-1c} tells us what the product along each column ought to be and~\ref{it:const-1d} tells us what the product along each row ought to be.

    A necessary and sufficient condition for compatibility of~\ref{it:const-1c} and~\ref{it:const-1d} is that the product of all the row products must equal the product of all the column products, as they are both the product of all entries in $M$.

    Let us now branch into three cases based on the additional conditions from~\cref{thm:setzer-crit} which are relevant to each: ${d \equiv \pm 1 \mod 8}$, ${d \equiv \pm 3 \mod 8}$, and $d \equiv \pm 2 \mod 8$. 

    \begin{enumerate}
        \item $\boldsymbol{d \equiv \pm 1 \mod 8}$: In this case, we have
            \begin{equation}
                \label{eq:rix-const--1}
                \prod_{q \mid n} \prod_{p \mid d} \leg{q}{p} = \prod_{q \mid n} \chi_d(q) = 1,
            \end{equation}
        where we have used that $\prod_{p \mid d} \leg{2}{p} = \boldsymbol{\chi}_{11}(|d|) = 1$. 

        We have incompatibility of~\ref{it:const-1c} and~\ref{it:const-1d} if 
        \[\prod_{p \mid d} \leg{-\sgn(d)\eps_d}{p} = -1,\]
        which happens if and only if $\sgn(d)\eps_d = 1$ and an odd number of the $p$ are congruent to $3 \mod 4$. The latter condition is equivalent to $\sgn(d)\eps_d = -1$, so this is never an issue.

        Because we are ordering the quadratic fields $\Q(\sqrt{m})$ by discriminant rather than by $m$, we must also keep track of the congruence class of $n$ modulo $4$ so that when $m = dn\sgn(d) \equiv 2,3 \mod 4$ we only take values of $n$ up to $\f{x}{4|d|}$. We thus distinguish the (sub)cases where $nd\sgn d \equiv 1 \mod 4$, $nd\sgn d \equiv 3 \mod 4$, and $n$ even. Denote by $R_{d}(x)\big|_{a (k)}$ the number of $n$ such that ${dn\sgn(d) \equiv a \mod k}$ and such that $\Q(\sqrt{dn\sgn(d)})$ has discriminant at most $x$ and admits an elliptic curve with GRE$_\Q$.
            Let us start with $R_d(x)\big|_{1 (4)}$, corresponding to the additional condition 
            \begin{enumerate}[(a), start = 5]
                \item $\leg{-4}{n} = \leg{-4}{d\sgn(d)}$.\label{it:const-1e1}
            \end{enumerate}

            Let $(b_n)_{n \geq 1}$ be the sequence of coefficients of the Dirichlet series
            \begin{equation}
                \label{eq:ridx-const-1}
                \qquad \f{1}{2}\paren{L(s,\boldsymbol{\chi}_2) + \leg{-4}{d\sgn(d)}L\paren{s,\leg{-4}{\argument}}} \, \otimes\, \bigotimes_{i=1}^{r-1} \frac{1}{2}\paren{\zeta(s) + \leg{-\sgn(d)\eps_d}{p_i}L\paren{s, \leg{\argument}{p_i}}},
            \end{equation}
            where we write $\boldsymbol{\chi}_2$ to denote the (principal) Dirichlet character with modulus $2$.  
            Note that the big convolution is over only $r-1$ primes -- we (arbitrarily) omit one of the prime factors of $d$. Consider any term in the big convolution. It corresponds to a Dirichlet series with a sequence of coefficients whose $n^{\t{th}}$ term -- assuming $n$ is coprime to $d$ -- is $1$ if $\leg{n}{p} = \leg{-\sgn(d)\eps_d}{p}$ and $0$ otherwise. Because $(a_n)_{n \geq 1}$ is $0$ whenever $n$ and $d$ are not coprime by~\ref{it:const-1a}, we do not need to worry about the behavior of $(b_n)_{n \geq 1}$ when $n$ and $d$ are not coprime.

            The term outside the big convolution corresponds to a Dirichlet series with a sequence of coefficients whose $n^{\t{th}}$ term is $1$ if \[\leg{-4}{n} = \leg{-4}{d\sgn(d)}\] and is $0$ otherwise. 

            Assuming that $n$ satisfies~\ref{it:const-1a}-\ref{it:const-1c}, we see that that $(b_n)_{n \geq 1}$ is constructed so as to be $1$ if and only if~\ref{it:const-1e1} is satisfied (because of the first term in~\eqref{eq:ridx-const-1}) and~\ref{it:const-1d} is satisfied for all but one of the primes dividing~$d$, and to be $0$ otherwise. However, we actually have more than this. Since 
            \[\prod_{p_i \mid d} \leg{-\sgn(d)\eps_d}{p_i} = 1,\]
            we see that knowing $\leg{n}{p_i}$ for $1 \leq i \leq r-1$ tells us $\leg{n}{p_r}$. We are exploiting here that we know that~\ref{it:const-1c} and~\ref{it:const-1d} are compatible. Thus, $a_n = b_n = 1$ if and only if~\ref{it:const-1a}-\ref{it:const-1e1} are satisfied. We have that
            \[R_d(x)\big|_{1 (4)} = \sum_{n \leq x} a_nb_n.\] 
            We will access this by applying~\cref{thm:selberg-delange} to 
            \[
            \sum_{n \geq 1} a_nb_nn^{-s} = F_d(s) \otimes \sum_{n \geq 1} b_nn^{-s}.
            \]
            Expanding out the convolutions of $\sum_n b_nn^{-s}$ in~\eqref{eq:ridx-const-1} yields a sum of Dirichlet $L$-series, each associated with a product of Kronecker characters. By~\cref{lem:ridx-factor-conv}, convolving $F_d(s)$ with $L(s, \chi)$ for any Dirichlet character $\chi$ gives
            \begin{equation}
                \label{eq:ridx-const-25}
                C(s)L(s, \chi)^{\nicefrac{1}{2}}L(s, \chi\chi_d)^{\nicefrac{1}{2}}
            \end{equation}
            where $C(s)$ is holomorphic on $\Re(s) > \frac{1}{2}$. \eqref{eq:ridx-const-25} will not have a singularity at $s = 1$ unless at least one of $\chi$ or $\chi\chi_d$ is principal. Because $\chi_d$ is primitive and real, $\chi\chi_d$ is principal if and only if $\chi$ is an extension by zero of $\chi_d$. This can happen only if the modulus of $\chi$ is divisible by $m_d$. However, none of the characters in~\eqref{eq:ridx-const-1} are zero at the omitted prime factor. Thus~\eqref{eq:ridx-const-25} has a singularity only when $\chi$ is principal. By~\cref{lem:ridx-factor-conv}, the term of interest is
                \[
                \begin{split}
                F_d(s) \otimes \frac{1}{2^{\omega(d)}}L(s,\boldsymbol{\chi}_{2}) =\ 
                & \frac{1}{2^{\omega(d)}}L(s,\boldsymbol{\chi}_2)^{\nf{1}{2}}
                L(s,\boldsymbol{\chi}_2\chi_d)^{\nf{1}{2}} L(2s,\boldsymbol{\chi}_2)^{-\nf{1}{2}} 
                \\ &
                \times \prod_{q | m_d} (1+\boldsymbol{\chi}_2(q)q^{-s}) \prod_{q:\chi_d(q)=1}(1-\chi^2(q)q^{-2s})^{\nf{1}{2}}.
                \end{split}
                \]

            Observing that 
            \begin{align*}
            L(s,\boldsymbol{\chi}_2)^{\nf{1}{2}} & = \zeta(s)^{\nf{1}{2}}(1-2^{-s})^{\nf{1}{2}}, \\
            L(s,\boldsymbol{\chi}_2\chi_d)^{\nf{1}{2}} & = L(s,\chi_d)^{\nf{1}{2}}(1-2^{-s})^{\nf{1}{2}},\\
            L(2s,\boldsymbol{\chi}_2)^{-\nf{1}{2}} & = \zeta(2s)^{-\nf{1}{2}}(1-2^{-2s})^{-\nf{1}{2}},
            \end{align*}
            $\boldsymbol{\chi}_2(q) = 1$ for all $q | m_d$, and $\chi_d(2) = 1$, we see that
            \[F_d(s) \otimes \frac{1}{2^{\omega(d)}}L(s,\boldsymbol{\chi}_{2}) = \f{1}{2^{\omega(d)}}(1+2^{-s})^{-1}F_d(s).\]
            Now, applying~\cref{cor:ridx-asymp} gives 

            \[R_d(x)\big|_{1 (4)} = \sum_{n \leq x} a_nb_n = \f{\f{2}{3}(1+o'_d(1))c_dx}{|d|2^{\omega(d)}\sqrt{\log x}}\]
            where the sum is up to $x$ because in this case $m = nd\sgn(d) \equiv 1 \mod 4$.

            We can run a similar argument for $R_d(x)|_{3 (4)}$. In this case, we can take $\sum_{n \geq 1} b_nn^{-s}$ to be
            \begin{align*}
                \qquad \f{1}{2}\paren{L(s,\boldsymbol{\chi}_2) - \leg{-4}{d\sgn(d)}L\paren{s,\leg{-4}{\argument}}} \, \otimes\, \bigotimes_{i=1}^{r-1} \f{1}{2}\paren{\zeta(s) + \leg{-\sgn(d)\eps_d}{p_i}L\paren{s, \leg{\argument}{p_i}}},
            \end{align*} 
            but since the only surviving term corresponds to $L(s,\boldsymbol{\chi}_2)$ the sign change in the first term of the product is irrelevant. The same argument as above yields
            \[R_d(x)\big|_{3 (4)}=\sum_{n \leq \f{x}{4}} a_nb_n = \f{\f{2}{3}(1+o'_d(1))c_dx}{4|d| 2^{\omega(d)}\sqrt{\log x}},\]
            where we only consider $n$ up to $\f{x}{4}$ because there is an additional factor of $4$ in $\Delta_{\Q(\sqrt{m})}$ when $m = dn\sgn(d) \equiv 3 \mod 4$.

            When computing $R_d(x)\big|_{2 (4)}$ we have the added condition    
            \begin{enumerate}[(a), start = 5]
                \item $n$ is even.
            \end{enumerate}

            We can filter out odd $n$ by convolving with $\zeta(s) - L(s,\chi_2)$, meaning that we can take $\sum_{n \geq 1} b_nn^{-s}$ to be
            \begin{align*}
                \qquad \big(\zeta(s) - L(s,\boldsymbol{\chi}_2)\big) \, \otimes \, \bigotimes_{i=1}^{r-1} \f{1}{2}\paren{\zeta(s) + \leg{-\sgn(d)\eps_d}{p_i}L\paren{s, \leg{\argument}{p_i}}},
            \end{align*}  As such, we have
            \[F_d(s) \otimes \f{1}{2^{\omega(d)-1}}\zeta(s) -F_d(s) \otimes L(s,\boldsymbol{\chi}_2).\]
            Thus, by~\cref{cor:ridx-asymp},
            \[R_d(x)\big|_{2 (4)} = \f{\f{1}{3}(1+o'_d(1))c_dx}{4|d| \cdot 2^{\omega(d)-1}\sqrt{\log x}}.\]

            Adding together the contributions from $1,2,$ and $3 \mod 4$, we see that 
            \[R_d(x) = \f{(1+o'_d(1))c_dx}{|d|2^{\omega(d)}\sqrt{\log x}}\]
            
        \item $\boldsymbol{d \equiv \pm 3 \mod 8}$: Conditions~\ref{it:const-1a}-\ref{it:const-1d} still apply. In addition, we require that 
            \begin{enumerate}[(a), start = 5]
                \item $m \equiv 1 \mod 4$.\label{it:const-1e3}
            \end{enumerate}
            Note that even though $\boldsymbol{\chi}_{11}(d) \neq 1$,~\ref{it:const-1e3} ensures that $n$ is odd and hence~\eqref{eq:rix-const--1} still holds, meaning that ~\ref{it:const-1a}-\ref{it:const-1d} are compatible. Additionally,~\ref{it:const-1e3} ensures that we need to consider just the case $R_d(x)\big|_{1(4)}$, and the computation is the same as in the $d \equiv \pm 1 \mod 8$ case, yielding

            \[R_d(x) = \f{\f{2}{3}(1+o'_d(1))c_dx}{|d|2^{\omega(d)}\sqrt{\log x}}.
            \]
        \item $\boldsymbol{d \equiv \pm 2 \mod 8}$: As usual,~\ref{it:const-1a}-\ref{it:const-1d} still apply. We also have the additional condition 
            \begin{enumerate}[(a), start = 5]
                \item $n \equiv d+1 \mod 8$.\label{it:const-1e2}
            \end{enumerate}

            We start with $d \equiv 6 \mod 8$ as it is simpler. We have by~\ref{it:const-1e2} that $\boldsymbol{\chi}_{11}(n) = 1$ in this case and hence~\eqref{eq:rix-const--1} holds. In addition, 
            \begin{equation}
                \label{eq:ridx-const-3}
                \prod_{i = 1}^r \leg{\sgn d}{p_i} = 1,
            \end{equation}
            because even if $\sgn d = -1$, that $d \equiv 6 \mod 8$ implies that we have $d' \equiv 3 \mod 4$ and hence an even number of the primes dividing $d'$ are $3 \mod 4$ (since $d < 0$). Thus, we have compatibility of~\ref{it:const-1c}-\ref{it:const-1e2}.

            As before, we now make the stronger claim that we only need to check~\ref{it:const-1c},~\ref{it:const-1e2}, and~\ref{it:const-1d} at all but one (odd) prime dividing $d$ to ensure that all three are satisfied. Suppose that $\chi_d(n) = 1$, $n \equiv 7 \mod 8$, and for $1 \leq i \leq r-1$,
            \[\leg{n}{p_i} = \leg{\sgn d}{p_i}.\]
            Then, we see that for such $n$ 
            \begin{align*}
                1 = \chi_d(n) 
                = \boldsymbol{\chi}_{11}(n) \prod_{1 \leq i \leq r} \leg{n}{p_i} 
                = \leg{n}{p_r}\prod_{i = 1}^{r-1} \leg{\sgn d}{p_i} 
            \end{align*}
            and comparing with~\eqref{eq:ridx-const-3} we have that 
            \[\leg{n}{p_r} = \leg{\sgn d}{p_r}.\]
            It is sufficient to use $\sum_n b_nn^{-s}$ to check that $n \equiv 7 \mod 8$ and the fourth condition for all but the last odd prime dividing $d$. We already know how to enforce the quadratic residuosity conditions. To force $n \equiv 7 \mod 8$, consider the function $\mathbf{1}_{7 (8)} \colon (\Z/8\Z)^{\times} \rightarrow \C$ that is $1$ on $7 \mod 8$ and $0$ everywhere else. The Fourier expansion of this function is
            \[\mathbf{1}_{7 (8)} = \tfrac{1}{4}(\boldsymbol{\chi}_{00} - \boldsymbol{\chi}_{01} - \boldsymbol{\chi}_{10} + \boldsymbol{\chi}_{11}).\]

            The initial factor we will add to our expression for $\sum_{n} b_nn^{-s}$ in this case is thus
            \[\tfrac{1}{4}(\zeta(s) - L(s, \boldsymbol{\chi}_{01}) - L(s, \boldsymbol{\chi}_{10}) + L(s, \boldsymbol{\chi}_{11})),\]
            where we can use $\zeta(s)$ instead of $L(s, \boldsymbol{\chi}_{00})$ as the coefficients of this factor at even indices are irrelevant. Therefore, we may take $\sum_{n \geq 1} b_nn^{-s}$ to be
            \[
            \qquad
            \tf{1}{4} \big(\zeta(s) - L(s, \boldsymbol{\chi}_{01}) - L(s, \boldsymbol{\chi}_{10}) + L(s, \boldsymbol{\chi}_{11})\big) \, \otimes \, \bigotimes_{i=1}^{r-1} \tf{1}{2}\paren{\zeta(s) + \leg{\sgn d}{p_i}L\paren{s, \leg{\argument}{p_i}}}.\]
            Following the same argument as before, we see that the only term of $\sum_n b_nn^{-s}$ that yields a term with a singularity at $s = 1$ after convolving with $\sum_n a_nn^{-s}$ is the term corresponding to $\zeta(s)$. Therefore, 
            \[R_d(x) = \sum_{n \leq \f{x}{4}} a_nb_n = \f{(1+o'_d(1))c_dx}{4|d|2^{\omega(d)}\sqrt{\log x}},\]
            noting in this case that $r := \omega(d) - 1$.

            When $d \equiv 2 \mod 8$ the situation is less simple. By~\cref{lem:chid} we have 
            
            \[1 = \chi_d(n) = \boldsymbol{\chi}_{01}(n) \prod_{1 \leq i \leq r} \leg{n}{p_i} = -\prod_{q | n}\prod_{p | d} \leg{q}{p},\]
            meaning that~\ref{it:const-1a}-\ref{it:const-1e2} are incompatible when 
            \[\prod_{i = 1}^r \leg{-\sgn(d)}{p_i} = 1.\]
            This happens when $\sgn(d) < 0$. If we assume $\sgn(d) > 0$ then we have compatibility and get the same answer as we got when $d \equiv 6 \mod 8$. 
    \end{enumerate}

    In the imaginary case,~\ref{it:const-1a}-\ref{it:const-1c} are the same but we instead have 
    \begin{enumerate}[(a), start = 4]
        \item \label{it:const-d-im} $\leg{n}{p} = \leg{\sgn(d)\eps_d}{p}$ for every odd $p \mid d$.
    \end{enumerate}

    The argument when $d$ is odd is identical except that the compatibility between~\ref{it:const-1c} and~\ref{it:const-d-im} now plays a role.~\ref{it:const-d-im} forces \[\leg{n}{p} = \leg{\sgn(d)\eps_d}{p}\] for every odd $p \mid d$. There is no way that~\ref{it:const-1c} can also be satisfied if $\sgn(d)\eps_d = -1$.\footnote{As an aside, notice that this is exactly the fifth constraint of~\cref{thm:setzer-crit}. This means that the fifth constraint of~\cref{thm:setzer-crit} is redundant for odd $d$.} When $d$ is even,~\ref{it:const-1c} and~\ref{it:const-d-im} are incompatible only when $d > 0$ and $d \equiv 6 \mod 8$. We also have the restriction from~\ref{thm:setzer-crit}\ref{it:setzer-5}, forbidding $d > 0$ when $d$ is even. 
\end{proof}

\begin{remark}
    In principle, we could have used sequences $(b_n)_{n \geq 1}$ which include all the prime factors of $d$. This approach yields multiple terms with singularities after convolution because our character expansion of $\sum_{n} b_nn^{-s}$ has terms which are Dirichlet $L$-series for characters $\chi$ induced by $\chi_d$. Our approach simplifies the computation. 
\end{remark}

\section{A Matching Lower Bound: Proving \texorpdfstring{\cref{thm:lcm-sparse}}{Theorem C} and \texorpdfstring{\cref{thm:rix-const}}{Theorem A}}
\label{sec:lb}

Consider the upper bound on $R(x)$ (the same argument works for $I(x)$ as well). Using~\cref{lem:ridx-const} instead of~\cref{lem:ridx-asymp} in the proof of~\cref{thm:rix-asymp} tells us that 
\begin{equation}
    \label{eq:lb-1}
    R(x) \leq \sum_{d \t{ good}} R_d(x) \sim \f{x}{\sqrt{\log x}} \sum_{d \t{ good}} \f{c_dc'_d}{|d|}.
\end{equation}

This is a priori only an upper bound because we may be double-counting -- if $d_1 \neq d_2$ are both good and both divide $m$ then $\Q(\sqrt{m})$ may be double-counted in~\eqref{eq:lb-1} since it could be that $\f{m}{d_1}$ contributes to $R_{d_1}(x)$  and $\f{m}{d_2}$ contributes to $R_{d_2}(x)$. Write $R_{d_1, d_2}(x)$ to denote the set of positive numbers $n$ such that $K = \Q(\sqrt{\lcm(d_1, d_2)n})$ admits an elliptic curve with good reduction everywhere and $\Delta_K$ is at most $x$. Then, by inclusion-exclusion, we have the lower bound
\begin{equation}
    \label{eq:lb-2}
    R(x) \geq \sum_{d \t{ good}} R_d(x) - \sum_{d, d' \t{ good}} R_{d,d'}(x).
\end{equation}

The idea is that this second term ends up being $\asymp \f{x}{\log^{\nicefrac{3}{4}} x}$ and hence is negligible compared to the first term (which is $\asymp \f{x}{\sqrt{\log x}}$ by~\cref{thm:rix-asymp}). As in the proof of~\cref{thm:rix-asymp}, there are two parts to this result. We first show that the sum of $\f{1}{\lcm(d,d')}$ over pairs of good $d$ and $d'$ converges. Then, we show that the dependence on $x$ of any $R_{d,d'}(x)$ is $\asymp \f{x}{\log^{\nicefrac{3}{4}} x}$. For the former, we show that the number of pairs $(d,d')$ where $d$ and $d'$ are good, have absolute value at most $x$, and $\lcm(d,d') \leq x$ is $x^{1-\kappa+o(1)}$ for some $\kappa > 0$. Then, the first part of~\cref{lem:nat-to-harm} implies that the sum of reciprocals of least common multiples with multiplicity is some absolute constant.

\lcmsparse

As motivation for~\crefrestated{thm:lcm-sparse}, notice that the analogous result with pairwise least common multiple replaced by pairwise product holds with $\f{\beta}{2-\beta}$ in place of $\beta$. To see this, consider splitting the range $[1,x]$ into intervals $(y, 2y]$. For each $a \in (y,2y]$ for which $a$ and $b$ are in $S$, any $b \in [1,x]$ such that $ab \leq x$ is at most $\f{x}{y}$. There are at most $O_{\veps}(x/y)^{1-\beta+\veps}$ such values of $b$ for any $\veps > 0$.
The number of possible values of $a \in (y,2y]$ is at most $O_{\veps}(y^{1-\beta+\veps})$. Therefore, the number of tuples $(a,b)$ where $a \in (y,2y]$ is at most \[O_{\veps}(y^{1-\beta+\veps}(x/y)^{1-\beta+\veps}) = x^{1-\beta+o(1)}\]
as $x$ goes to infinity. This is uniform in $y$ and there are at most $\log x + 1$ intervals $(y,2y]$. Therefore, the total number of pairs $(a,b)$ which work is at most $x^{1-\beta+o(1)}(\log x + 1) = x^{1-\beta+o(1)}$. The structure of our proof of~\crefrestated{thm:lcm-sparse} is similar. The main difficulty is that $\lcm(a,b)$ may be much smaller than $ab$ so it is harder to control the number of $b$ which can be associated to a given $a$.

We present an improved proof of~\crefrestated{thm:lcm-sparse} due to Ashwin Sah and Mehtaab Sawhney, and we thank them for allowing us to present it in this paper.\footnote{The original proof of the result, due to the authors, gave only an upper bound and had a slightly weaker exponent.} We will use the following standard lemma. 

\begin{lemma}
    \label{lem:divisor-bound}
The number of divisors of $n$ is $\exp(O(\log n/\log\log n)) = n^{o(1)}$.
\end{lemma}

\begin{proof}[Proof of~\crefrestated{thm:lcm-sparse}]
    Let $\kappa := \beta/(2-\beta)$. We wish to bound the total number of pairs $(a,b) \in S\times S$ with $\lcm(a,b)\le x$. We consider such pairs with $a\in[y,2y), b\in[z,2z)$, $\gcd(a,b) \in[g,2g)$ for some $yz \geq x$ (if $yz \leq x$ then we can just use the product argument above). Note that $g\le\min(2y,2z)$. We will show that the number of pairs with $a,b,$ and $g$ in these ranges is at most $x^{1-\kappa+o(1)}$. Then, summing the contributions from all such triples of intervals only adds a factor of $O(\log^3 x) = x^{o(1)}$ to the overall bound if we take a dyadic decomposition.

    First, observe that
    \[x\ge\lcm(a,b) = \frac{ab}{\gcd(a,b)}\ge\frac{yz}{2g}.\]

    For any $\veps > 0$, there are at most $O_{\veps}(y^{1-\beta+\veps})$ choices of $a$ and at most $O_{\veps}(z^{1-\beta+\veps})$ choices of $b$ by $\beta$-polynomial sparsity of $S$. Therefore, there are at most $O_{\veps}((yz)^{1-\beta+\veps})$ choices of pairs. Because $yz \leq 2xg$, this is 
    \begin{equation}
        \label{eq:lcm-sparse-1}
        O_{\veps}((xg)^{1-\beta+\veps}).
    \end{equation}

    We can bound the number of pairs another way. There are at most $O_{\veps}(y^{1-\beta+\veps})$ choices of $a$. For each such $a$, there are then $O_{\veps}(x^{\veps})$ divisors of $a$ lying in $[g,2g]$ (i.e.\ choices for $\gcd(a,b)$) by~\cref{lem:divisor-bound}. There are $O(z/g)$ choices of $b$ divisible by this choice of $\gcd(a,b)$. This gives a bound of $O_{\veps}(x^{\veps}y^{1-\beta}z/g)$.

    We may obtain a symmetric bound by swapping the roles of $a,b$, giving a bound of
    \[x^{o(1)}\cdot \min\paren{y^{1-\beta}z/g,yz^{1-\beta}/g}\]
    as $x$ goes to infinity. Taking the geometric mean of the two terms in the minima gives
    \begin{equation}
        \label{eq:lcm-sparse-2}
        x^{o(1)}(yz)^{1-\nf{\beta}{2}}/g = x^{1-\nf{\beta}{2}+o(1)}/g^{\nf{\beta}{2}},
    \end{equation}
    using again that $yz\le 2xg$.

    Combining~\eqref{eq:lcm-sparse-1} and~\eqref{eq:lcm-sparse-2}, we obtain a bound of 
    \[\min\big(O_{\veps}((xg)^{1-\beta+\veps}),O_{\veps}(x^{1-\nf{\beta}{2}+\veps})/g^{\nf{\beta}{2}}\big).\]
    This is maximized when the two terms are approximately equal, which happens when $g$ becomes $x^{\beta/(2-\beta)+o(1)}$ as $x$ goes to infinity. This yields the bound
    \[x^{1-\f{\beta}{2-\beta}},\]
    implying $\beta/(2-\beta)$-polynomial sparsity.

    For the matching lower bound, consider the set $S$ constructed as follows. Pick some positive integer $x_0$ and add to $S$ the multiples of $\ceil{x_0^{\kappa}}$ in the interval $\big[\f{1}{2}x_0^{(1+\kappa)/2}, x_0^{(1+\kappa)/2}\big)$. Then, we have added at most $x_0^{(1-\kappa)/2}$ values to $S$, the least common multiple of any pair of such values is at most $x_0$, and the number of tuples of elements is $x_0^{1-\kappa}$. Continue by choosing $x_1$ much larger than $x_0$ and repeating the process for each $x_i$ for all $i \in \NN$.

    The number of elements of $S$ up to $x$ grows as 
    \[x^{\f{1-\kappa}{1+\kappa}} = x^{1-\beta},\]
    and the number of least common multiples up to $x$ grows as $x^{1-\kappa+o(1)}$, so we see that we have a lower bound matching our upper bound.  
\end{proof}

Applied to the set of good $d$, which by~\crefrestated{thm:red-ec} satisfies the conditions of~\crefrestated{thm:lcm-sparse} for $\beta = \nf{1}{3}$, we can take $\kappa = \nf{1}{5}$.  

\begin{lemma}
    \[R_{d,d'}(x) \ll \f{(1+o'_{dd'}(1))c_{dd'}x}{\lcm(|d|,\!|d'|)\log^{\nf{3}{4}} x} \quad\t{and}\quad I_{d,d'}(x) \ll \f{(1+o'_{dd'}(1))c_{dd'}x}{\lcm(|d|,\!|d'|)\log^{\nf{3}{4}} x},\]
    where $c_{dd'} \ll \lcm(|d|,|d'|)^{0.001}$ and $o'_{dd'}(1)$ denotes some function of $x$ and $d$ which goes to $0$ as $x$ goes to infinity while $\lcm(d,d) \leq \log^{k}x$ for some constant $k$. 
\end{lemma}
\begin{proof}
    As in the proof of~\cref{lem:ridx-asymp}, we will use only the first constraint from~\cref{thm:setzer-crit} and order quadratic fields $\Q(\sqrt{m})$ by $m$ rather than by discriminant. Let $S$ (resp.\ $S'$) be the set of primes $q$ such that $\chi_d(q) = 1$ (resp.\ $\chi_{d'}(q) = 1$). An upper bound on $R_{d, d'}(x)$ is the number of squarefree $n$ at most $\f{x}{\lcm(|d|,|d'|)}$ which are divisible only by primes in $S \cap S'$.\footnote{Note that this condition is necessary but not sufficient, as $d$ and $d'$ must themselves be compatible in the sense that the primes dividing $\f{d'}{\lcm(|d|,|d'|)}$ must lie in $S$ and vice-versa. However, even this weaker condition is enough.} We have that when $q$ is coprime to $d$ and $d'$, 
    \[\f{1}{4}(\chi_d(q) + 1)(\chi_{d'}(q) + 1) = \begin{cases}
        1 & q \in S \cap S' \\
        0 & q \notin S \cap S' \\
    \end{cases}
    \]

    The expression on the left-hand-side can be written as
    \[\f{1}{4}(1 + \chi_d(q) + \chi_{d'}(q) + \chi_d\chi_{d'}(q)).\]
    Let 
    \[a_n := \begin{cases}
        1 & n \t{ is squarefree and divisible only by primes in $S \cap S'$} \\
        0 & \t{otherwise.} \\
    \end{cases}\]

    Note that $\chi_d$ and $\chi_{d'}$ are primitive (or primitive up to a local factor, when $d \equiv \pm 3 \mod 8$) and nonprincipal. We have as in the proof of~\cref{lem:ridx-asymp} that 
    \[F(s) := \sum_{n \geq 1} a_nn^{-s} = C_{d,d'}(s)\big(\zeta(s)L(s,\chi_d)L(s,\chi_{d'})L(s,\chi_{d}\chi_{d'})\big)^{\nf{1}{4}},\]
    where $C_{d, d'}(s)$ is holomorphic on $\Re(s) > \f{1}{2}$ and $C_{d,d'}(1) \ll \lcm(|d|,|d'|)^{0.001}$. Because the convolution of primitive quadratic characters is principal if and only if they are equal, we see that we can apply~\cref{thm:selberg-delange} with 
    \[F(s) = \zeta^{\nf{1}{4}}(s)G(s)\]
    for 
    \[G(s) = \big(L(s,\chi_d)L(s,\chi_{d'})L(s,\chi_{d}\chi_{d'})\big)^{\nf{1}{4}}C_{d,d'}(s).\]
    This gives us the desired bound on the sum of $a_n$ up to $\f{x}{\lcm(|d|,|d'|)}$ after some manipulation,~\cref{cor:robin}, and~\cref{lem:pl-l-func} as before. 
\end{proof}

\begin{lemma}
     \[\sum_{d, d' \textnormal{ good}} R_{d,d'}(x) \ll \f{x}{\log^{\nf{3}{4}} x} \] 
\end{lemma}

\begin{proof}
    This follows from the same argument as was used in deriving~\cref{thm:rix-asymp} from~\cref{lem:ridx-asymp}. We have for any $z \leq x$ that 
    \begin{align*}
        \sum_{d, d' \t{ good}} R_{d,d'}(x) & \ll x\paren{\sum_{\substack{d, d' \t{ good} \\ \lcm(|d|,|d'|) \leq z}} \f{(1+o'_{dd'}(1))c_{dd'}}{\log^{\nf{3}{4}} x} + \sum_{\substack{d,d' \t{ good} \\ \lcm(|d|,|d'|) > z }} \f{1}{\lcm(|d|,|d'|)}}. \\
    \end{align*}
    Taking $z := \log^4 x$, we see that the second term in the parenthesis is negligible compared to the first as $x$ goes to infinity by~\crefrestated{thm:red-ec}, \crefrestated{thm:lcm-sparse} and~\cref{lem:nat-to-harm}\ref{it:nat-to-harm-2}. Every $o'_{dd'}(1)$ for $|d| \leq z$ is then upper bounded by some $o(1)$ independent of $d$ and $d'$. Therefore, as $x$ goes to infinity we have an asymptotic bound 
    \[\ll \f{x}{\log^{\nf{3}{4}} x}\sum_{\substack{d,d' \t{ good}}} \f{c_{dd'}}{\lcm(|d|,|d'|)}\ll \f{x}{\log^{\nf{3}{4}} x}\]
    as desired.
\end{proof}

\cref{thm:rix-const} then follows from~\cref{lem:ridx-const}, \eqref{eq:lb-1} and \eqref{eq:lb-2}, and~\cref{thm:lcm-sparse}.

\section{Computing the Constants: \texorpdfstring{Proving~\cref{cor:const-lb} and~\cref{cor:const-ub}}{Proving Corollary D and Corollary E}}
\label{sec:comp}

We now turn to the problem of obtaining numerical estimates for $c_R$ and $c_I$ in~\cref{thm:rix-const}, or equivalently, for series of the form 
\[\sum_{d \t{ good}} \f{c_dc'_d}{|d|2^{\omega(d)}}.\]
where $c_d'$ is one of $c'_{d,R}$ or $c'_{d,I}$.

Our approach is to first compute this series explicitly for good $d$ up to $\abs{d} \leq D$ and to then bound the size of the tail. This would give us an estimate along with an error bound. However, while~\crefrestated{thm:red-ec} does tell us that the tail has size $\ll_{\eps} D^{-\nf{1}{3}+o(1)}$ as $D$ goes to infinity, our dependence on the $abc$-conjecture means that we cannot control the leading constant. Because each term is nonnegative, we can lower bound the sum by ignoring the tail, yielding~\cref{cor:const-lb}. However, we need some more information in order to control the tail and show an upper bound.

\begin{remark}[Explicit $abc$-conjectures]
One might wonder if we can obtain the desired bounds via an explicit formulation of the $abc$-conjecture. For example, Robert, Stewart, and Tenenbaum~\cite{explicit-merit-abc} conjectured that
\begin{equation}
    \max(|a|,|b|,|c|) < k^{1+\eps(k)},
\end{equation}
where $k$ is the radical of $a,b,c$ in~\cref{conj:abc},
\[\veps(k) := \sqrt{\f{48}{\log k\log\log k}}\paren{1+\f{3\log\log\log k+2C_1}{2\log\log k}},\]
and $C_1 := 1+\log 3 - \f{13}{6}\log 2 + \veps$ for any $\veps > 0$.

Consider $(r,d,t)$ such that $r^3 = dt^2 - 1728$, and let $k := \rad(1728r^3dt^2) = \rad(6rdt)$. Then, the explicit $abc$-conjecture applied to $(r^3, -dt^2, 1728)$ tells us that 
\[|r|^{\nf{1}{2}-\nf{5}{2}\veps(k)} \leq 6^{1+\veps(k)}|d|^{\nf{1}{2}+\nf{1}{2}\veps(k)}.\]

This is trivial unless $\veps(k) < \f{1}{5}$, which does not happen until $k > 10^{141}$. As such, we cannot hope for any useful bound on the tail until $k > 10^{141}$. The best lower bound we can presently prove on $k$ in terms of $d$ is that $k = \rad(6drt) \gg d$. Without a better lower bound, this seems to require that we explicitly compute the contributions to the constant for $d < 10^{141}$. We encounter similar obstacles when attempting to use Baker's explicit $abc$-conjecture~\cite{baker}; see for example the table in Theorem 1 of \cite{laishram-shorey}.
\end{remark}

\subsection{The frequency of good \texorpdfstring{$d$}{d}}

We start by motivating the assumption under which we will prove our upper bound. 
Recall Granville's conjecture on the twists of hyperelliptic curves (\cref{conj:granville-1}). Intuitively, it suggests that $\sqf(f(x))$ is usually not much smaller than $f(x)$. Here is the precise statement of the conjecture in the setting of elliptic curves. 

\begin{conjecture}[Granville~\cite{granville-twists}]
    \label{conj:granville-3}
    Let $E$ be an elliptic curve given by the integral model $y^2 = f(x)$, and write $f_3$ for the leading coefficient of $f$. Then, 
    \begin{equation}
        T_E(D) \sim \kappa_fD^{\nf{1}{3}},\label{eq:granville-3-1}
    \end{equation}
    where 
    \begin{equation}
        \label{eq:granville-3-2}
        \kappa_f := 2|f_3|^{-\nf{1}{3}}\prod_p \paren{1+\paren{1-\f{1}{p^{\nf{2}{3}}}}\paren{\f{\omega_f(p^2)}{p^{4/3}} + \f{\omega_f(p^4)}{p^{8/3}} + \f{\omega_f(p^6)}{p^{4}} + \dots}}
    \end{equation}
    and $\omega_f(r)$ is the number of roots of $f$ in $\Z/r\Z$. 
\end{conjecture}

Granville~\cite{granville-abc} also showed that the lower bounds implicit in~\cref{conj:granville-1} hold under the $abc$-conjecture (with the specified constants $\kappa_f$). Granville~\cite{granville-twists} also proved~\cref{conj:granville-1} for high genus hyperelliptic curves which split into linear factors. Implicit in this latter proof was the following. 

\begin{lemma}[Implicit in Granville~\cite{granville-twists}]
    Let $f(x) \in \Z[x]$ be a separable cubic polynomial. Assume that the number of $r$ for which $|\sqf(f(r))| \leq D$ and $f(r)/\sqf(f(r)) \geq U(D)^2$ is $o(D^{\nf{1}{3}})$ for some $U(D) \ll o(D^{\nf{1}{6}})$. Then,~\cref{conj:granville-3}\footnote{The same condition extended to higher degree polynomials implies~\cref{conj:granville-1} by a similar argument.} holds. 
\end{lemma}

We prove the following. 

\begin{theorem}
    \label{thm:granchild}
    Let $f(x) := x^3-1728$.
    Assume that the number of $r$ for which $|\sqf(f(r))| \leq D$ and $f(r)/\sqf(f(r)) \geq U(D)$ is $o(D^{\nf{1}{3}})$ for some $U(D) \ll o(D^{\nf{1}{6}})$.
    Let $G(D)$ denote the number of good $d$ for which $|d| \leq D$.
    Then, $G(D) \sim \kappa'D^{\nf{1}{3}}$ for an absolute constant $\kappa$ such that $3.48523 \leq \kappa \leq 3.50692$.
\end{theorem}

\begin{proof}
    As the proof largely follows the proof of Theorem $2$ in~\cite{granville-twists}, we will highlight the differences and leave some details to the reader. Let $E$ be the elliptic curve given by integral model $y^2 = f(x)$. The proof of Theorem $2$ shows, under the assumption of the theorem statement, that
    \begin{equation}
        \label{eq:granchild-1}
        |T_E(D)| = \sum_{t \leq U(D)} \#\{r \colon |f(r)| \leq Dt^2, t^2 | f(r), f(r)/t^2 \t{ is squarefree}\} + o(D^{\nf{1}{3}}),
    \end{equation}
    and proceeds to argue that the first term is $\kappa_fD^{\nf{1}{3}}$ as $D$ goes to infinity. Henceforth, let $E$ be defined by the integral model $y^2 = f(x) := x^3-1728$, and write $\cR$ to denote the set in~\eqref{eq:good}. We wish to evaluate the sum in~\eqref{eq:granchild-1} while restricting to those $r$ satisfying~\eqref{eq:good}. Granville's proof, which has no such restriction, breaks each term into a sum over $r_0$ for which $f(r_0) \equiv 0 \mod{t^2}$, and, defining \[g(s) = f(r_0 + t^2s)/t^2,\] counts the number of $s$ for which $g(s)$ is squarefree. As $D$ (and hence the range of valid $s$) goes to infinity, this number is asymptotically 
    \[\prod_p \paren{1-\f{\omega(g(x),p^2)}{p^2}},\]
    where $\omega(h(x), m)$ is the  number of congruence classes $j \mod m$ for which $h(j) \equiv 0 \mod m$.\footnote{Granville denotes this number by $\omega_h(m)$} This asymptotic is uniform across terms of the sum in~\eqref{eq:granchild-1} because $U(D)$ grows slowly as a function of $D$. A Chinese remainder theorem argument then lets us write~\eqref{eq:granchild-1} as a product of form
    \begin{equation}
        \label{eq:granchild-2}
        \prod_p \sum_{r_p:f(r_p) \equiv 0\,\mathrm{mod}\,p^2} \paren{1-\f{\omega(g(r),p^2)}{p^2}},
    \end{equation}
    where the function $g$ depends on $r_p$. Studying $\omega(g,p^2)$ simplifies the product, and some manipulation then yields~\eqref{eq:granville-3-2}. Much of Granville's proof works for our application: introducing constraints which restrict the $r$ to certain congruence classes modulo $2$ and $3$ changes only the local factors of the product in~\eqref{eq:granville-3-2} corresponding to $2$ and $3$. These factors contribute to those $t$ which are divisible by $2$ or $3$. We will compute the contribution from each possible value of $|t|_2$ and $|t|_3$ and aggregate them to obtain replacements for the local factors at $2$ and $3$ in~\eqref{eq:granville-3-2}.

    If $t$ is divisible by $2$ then certainly $r$ is even and therefore $f(r) = (16u+v)^3-1728$ for some $u \in \Z$ and $v \in \{0,4\}$. Then, for $v' \in \{0,1\}$, we have 
    \begin{align*}
        dt^2 & = f(r) = 64((4u+v')^3-27) = 64(4u+v'-3)(16u^2+8uv'+v'^2+12u+3v'+9).
    \end{align*}
    Regardless of the value of $v'$, we see that $|t|_2 = \f{1}{8}$ if $t$ is even and its contribution to the constant is nonzero. In this case, the local factor at $2$ in~\eqref{eq:granchild-2} is
    \begin{equation}
        \label{eq:granchild-3}
        \sum_{\substack{r_2 \colon f(r_2) \equiv \,0\,\mathrm{mod}\,64 \\ r_2 \in \cR}} \paren{1-\f{\omega(f(r_2+64s)/64,4)}{4}},
    \end{equation}
    where we say that a congruence class is in $\cR$ if every integer in this congruence class is in $\cR$.
    As in Granville's proof, the numerator in each term of~\eqref{eq:granchild-3} is the number of roots of $g(s) = f(r_2+64s) \mod{256}$ which are congruent to $r_2$ modulo $64$.

    Over $\Z/2^k\Z$ for $k \leq 6$, $\bar{f}(x) = x^3$ has $2^{k-\ceil{\f{k}{3}}}$ roots, corresponding to multiples of $2^{\ceil{\f{k}{3}}}$. For $k > 6$, $12$ is always a root of $\bar{f}$. This root is simple because $f'(12) = 432 \not\equiv 0 \mod{2^7}$. The polynomial $x^2+12x+144$ has no roots modulo $32$, let alone higher powers, as \[n^2 = 32k - 108 = 4(8k-27) \equiv 4(8k+5)\]
    and $5$ is not a square modulo $8$.

    The values of $r_2$ in~\eqref{eq:granchild-3} are $0$ and $2$. No $r = 2 + 64s$ can lie in $\cR$, so we take $r_2 = 0$. The relevant roots of $g(s) \mod{64}$ are the multiples of $4$ which are congruent to $0$ or $4 \mod{16}$. None of these lift to $12$, which is the unique root of $g(s) \mod{256}$. Thus, the value of~\eqref{eq:granchild-3} is $1$.

    If $t$ is not divisible by $2$ then the local factor associated with $2$ in~\eqref{eq:granchild-2} would be
    \[\sum_{\substack{r_2 \colon f(r_2) \equiv\,0\,\mathrm{mod}\,4}} \paren{1-\f{\omega(f(r_2+4s),4)}{4}}.\]
    We need to force $r = r_2 + 4s \in \cR$. As before, we only need to consider $r_2 = 0$. In this case, we only allow $s \equiv 0,1 \mod 4$ and so our local factor in~\eqref{eq:granchild-2} is $\f{1}{2}$.

    We study the local factor at $3$ in a similar fashion. If $|t|_3 \neq 1$, $r \equiv 12 \mod{27}$ and
    \begin{align*}
        dt^2 & = f(r) = 27((9u+4)^3-64) = 3^6u(27u^2+216u+1456).
    \end{align*}
    In this case, we see that $|t|_3 \leq 3^{-3}$. Say $|t|_3 = 3^{-e_3}$ for $e_3 \geq 3$. Then, 
    \begin{equation}
        \label{eq:granchild-4}
        \sum_{\substack{r_3 \colon f(r_3) \equiv\,0\,\mathrm{mod}\,3^{2e_3} \\ r_3 \in \cR}} \paren{1-\f{\omega(f(r_3+3^{2e_3}s)/3^{2e_3},9)}{9}}.
    \end{equation}

    Over $\Z/3^k\Z$ for $k \leq 3$, $\bar{f}(x) = x^3$ has $3^{k-\ceil{\f{k}{3}}}$ roots. For $k > 3$, the root at $12$ is simple, and $x^2+12x+144$ has no roots modulo $81$. For all roots of $f(x) \mod 3^{e_3}$ in $\cR$, there is always exactly one lift to a root of $\Z/3^{e_3+2}\Z$ because the only root is $12$. Thus, our local factor in~\eqref{eq:granchild-4} in all of these cases is $\f{8}{9}$. A similar argument shows that the local factor in~\eqref{eq:granchild-4} when $|t|_3 = 1$ is $\f{8}{9}$.

    Repeating the manipulations in Granville's proof tells us that the local factor at $2$ in~\eqref{eq:granville-3-2} is 
    \[\f{1}{2} \cdot 1 + 1 \cdot (2^{-6})^{\nf{1}{3}}  = \f{3}{4}\] 
    and the local factor at $3$ is 
    \[\f{8}{9}(1+3^{-\nf{6}{3}} + 3^{-\nf{8}{3}} + \dots) = \f{8}{9}\paren{1+\f{1}{9}(1-3^{-\nf{2}{3}})^{-1}}.\]

    It remains to study the local factors in~\eqref{eq:granville-3-2} at those primes which are not $2$ or $3$. By Hensel's lemma, $\omega_f(p) = \omega_f(p^k)$ for each $k \geq 1$ if $p$ does not divide the discriminant of $f$. Therefore, the local factor in~\eqref{eq:granville-3-2} for $p \neq 2, 3$ is
    \[1+\omega(f,p)\f{p^{\nf{2}{3}}-1}{p^2-p^{\nf{2}{3}}}.\]
    Since $f(x) = (x-12)(x^2+12x+144)$, $12$ is always a simple root of $\bar{f}$ modulo such primes $p$. The reduction of $x^2-12x+144$ splits if and only if $-108$ is a quadratic residue, or equivalently $-3$ is a quadratic residue. Thus, the local factor for $p \neq 2,3$ is 
    \[1+\paren{2+\leg{-3}{p}}\f{p^{\nf{2}{3}}-1}{p^2-p^{\nf{2}{3}}}.\]

    Multiplying these together for the primes less than $P := 10000$ (not including the aforementioned factors at $2$ and $3$) yields $\approx 2.1532$. The product from the remaining primes is at most 
    \[\prod_{p \geq P} \paren{1+\f{3}{p^{\nf{4}{3}}-1}} \leq \exp\bigg(3\sum_{p \geq P} (p-1)^{-\nf{4}{3}}\bigg) \leq \exp\paren{\f{3\pi(10000)}{9999^{\nf{4}{3}}} + 4\int_{10000}^{\infty} \f{\pi(x)}{(x-1)^{\nf{7}{3}}} dx} \leq 1.0063\]
    where we have used that $(P-1)^{\nf{4}{3}} \leq P^{\nf{4}{3}} - 1$ for $P \geq 1$. 
\end{proof}

In order to get useful upper bounds on $c_R$ and $c_I$ in~\cref{thm:rix-const} (and hence show that $c_R > c_I$), we do not need an asymptotic result, but instead an upper bound on the frequency of good $d$ for which $|d| \leq D$. We are unable to prove such a bound unconditionally (or using the $abc$-conjecture), but prove the result conditional on the the following assumption. 
\begin{assumption}
    \label{assumption} $G(D) \leq 5D^{0.35}.$
\end{assumption}
\cref{thm:granchild} tells us that $G(D) \sim \kappa'D^{\nf{1}{3}}$, and hence we expect in the limit that the coefficient is less than $5$ and the exponent is less than $0.35$. Per~\cref{fig:d-conj}, the experimental frequency of good $d$ appears to converge rapidly to the prediction of~\cref{thm:granchild}, and the upper bound given by~\cref{assumption} appears to hold quite comfortably. 


\definecolor{colorUpperBound}{rgb}{0.6,0.3,1.0}
\definecolor{colorScaledGranville}{rgb}{0.5,0.5,0.5}
\definecolor{colorScaledGranvilleFill}{rgb}{0.3,0.3,0.3}
\definecolor{colorData}{rgb}{0,0,1.0}

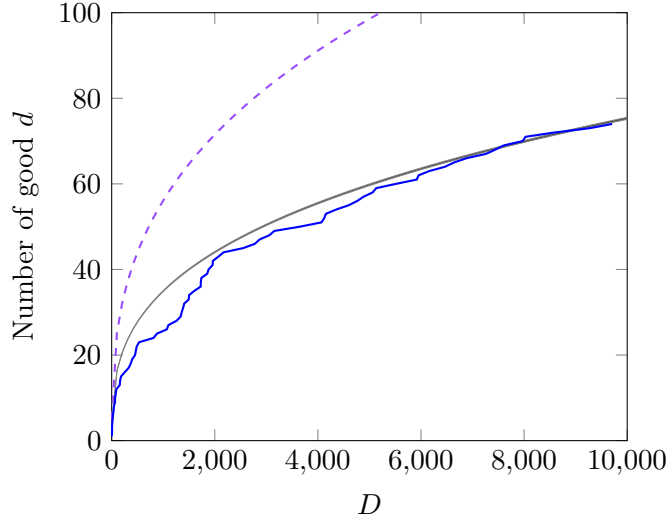
\begin{figure}[H]
    \pgfplotsset{scaled x ticks=false}
    \begin{tikzpicture}
        \begin{axis}[
            xlabel={$D$},
            xmin=0, xmax=10000,
            ylabel={Number of good $d$},
            ymin=0, ymax=100,
        ]
            \addplot[
                domain=1:10000,
                samples=100,
                color=colorUpperBound,
                thick,
                dashed,
                ]{5*x^(0.35)};
            \addplot[
                name path=A,
                domain=1:10000,
                samples=100,
                thin,
                color=colorScaledGranville,
                ]{3.48523*x^(1/3)};
            \addplot[
                name path=B,
                domain=1:10000,
                samples=100,
                thin,
                color=colorScaledGranville,
                ]{3.50692*x^(1/3)};
            \addplot[color=colorScaledGranvilleFill] fill between[of=A and B];
            \addplot[
                color=colorData, 
                thick,
            ] table [x=x,y=y, col sep=comma] {good-d-count.csv};
        \end{axis}
    \end{tikzpicture}
    \caption{The number of good $d$ with absolute value at most $D$ (solid blue), the bounds from the scaled version of Granville's conjecture in~\cref{thm:granchild} (solid grey), and our hypothesized upper bound (dashed purple).}
    \label{fig:d-conj}
\end{figure}

\subsection{An outline of the algorithm}
\label{subsec:algo}

We briefly outline the steps we use to compute our upper and lower bounds on $c_R$ and $c_I$. We direct the interested reader to the \href{https://github.com/abhijit-mudigonda/everywhere-good-reduction}{\color{blue} source code}.\footnote{\url{https://github.com/abhijit-mudigonda/everywhere-good-reduction}}

\subsubsection{Computing a list of good \texorpdfstring{$d$}{d}}
\label{subsubsec:good-d}

We compute a list of all good $d$ with $|d| \leq D$ by checking for each squarefree $d$, $|d| \leq D$, whether $E_d: dy^2 = x^3 - 1728$ has an integral point with $y\neq 0$ and $x$ an element of~\eqref{eq:good}.
To compute $E_d(\ZZ)$ we use the algorithm from~\cite{vonkaenel-matschke}, implemented in Sage~\cite{sagemath}, which is based on modularity and an elliptic logarithm sieve. 
It requires the knowledge of a Mordell--Weil basis for~$E_d$.
In most of our cases, the latter becomes the bottleneck of the computation of $E_d(\ZZ)$.
Note that $E_d$ is isomorphic over~$\Q$ to the Mordell curve $y^2=x^3 - 27d^3$, which in turn is $3$-isogenous to $E'_d: y^2 = x^3 + d^3$.
Thus it is enough to compute the Mordell--Weil basis for $E'_d$, to push it forward to $E_d$ via the $3$-isogeny, and to saturate it.

In most cases we use Magma~\cite{magma} to compute the Mordell--Weil basis for~$E_d$ directly. 
In case the rank of $E_d$ is $1$, it can be advantageous to find a Mordell--Weil basis with the Heegner point method, e.g.\ for~$d=131$:
In that case, first we estimate the regulators of $E_d$ and $E'_d$ (where the estimate is based on BSD), to see which of the two curves is expected to have a Mordell--Weil generator of smaller N\'eron--Tate height.
For that curve we compute a Heegner point using Pari/GP~\cite{pari}. 

\subsubsection{Computing the constant for each \texorpdfstring{$d$}{d}}
\label{subsubsec:each-cd}

For each such $d$, we compute 

\[\f{c_{d}c'_d}{|d|2^{\omega(d)}},\] 
where $c'_d$ is one of $c'_{d,R}$ or $c'_{d,I}$.

Given a value of $d$, the values of $c'_{d,R}, c'_{d,I}, |d|$, and $2^{\omega(d)}$ are easy to compute. Recall from~\cref{cor:ridx-factor-sqf} that  

\[c_d := \f{1}{\Gamma(\nicefrac{1}{2})\zeta(2)^{\nicefrac{1}{2}}}L(1,\chi_d)^{\nicefrac{1}{2}}\prod_{q \mid m_d} \paren{1+q^{-1}}^{-\nicefrac{1}{2}}\prod_{q:\chi_d(q) = 1} \paren{1-q^{-2}}^{\nicefrac{1}{2}}.\]
  
We compute $L(1,\chi_d)$ using the following theorem.

\begin{theorem}[Landau~\cite{davenport}]
Let $\Delta$ be a fundamental discriminant. Then, 

\[L(1, \leg{\Delta}{\argument}) = \begin{cases}
    -\f{\pi}{|\Delta|^{\nf{3}{2}}}\sum_{j = 1}^{|\Delta|} j\leg{\Delta}{j} & \Delta < 0 \\
    -\f{1}{|\Delta|^{\nf{1}{2}}}\sum_{j = 1}^{\Delta} j\log\sin\f{j\pi}{\Delta} & \Delta > 0 \\
    \end{cases}
\]
\end{theorem}

Because $L(s,\chi_d)$ agrees with a Kronecker character up to a constant factor (\cref{lem:chid}), this allows us to evaluate the L-series at $1$. 

To bound the product, we observe that for any $Q > 0$, 
\begin{equation}
    \zeta(2)^{-\nf{1}{2}}\prod_{\substack{q:\chi_d(q) \neq 1 \\ q \leq Q}} \paren{1-q^{-2}}^{-\nf{1}{2}} \leq \prod_{q:\chi_d(q) = 1} \paren{1-q^{-2}}^{\nicefrac{1}{2}} \leq \prod_{\substack{q:\chi_d(q) = 1 \\ q \leq Q}} \paren{1-q^{-2}}^{\nicefrac{1}{2}},
\end{equation}
and by taking $Q$ large enough we can obtain a decent approximation to this product. We could likely improve upon this by using the results of~\cite{explicit-dirichlet} but this is already enough to prove~\cref{cor:const-ub}.

\subsubsection{Bounding the size of the tail}
\label{subsubsec:tail}

To bound the contributions to the constants in~\cref{thm:rix-const} of good $d$ outside a given range, we first note that

\[\f{c_dc'_d}{|d|2^{\omega(d)}} \leq \f{L(1,\chi_d)^{\nf{1}{2}}}{2\Gamma(\nf{1}{2})\zeta(2)^{\nf{1}{2}}|d|}.\]

An upper bound on $L(1,\chi_d)$ comes from the following. 

\begin{lemma}[P\'olya-Vinogradov inequality~\cite{polya, davenport}]
    \label{lem:pv-ineq}
    Let $M$ and $N$ be positive integers. If $\chi$ is a primitive character with modulus $m$,
    \[\abs{\sum_{n = M}^{M+N} \chi(n)} < m^{\nf{1}{2}}\log m.\] 
\end{lemma}
\begin{corollary}
    \label{cor:pv-ineq}
    \begin{equation}
        \label{eq:pv-ineq}
        L(1,\chi_d) \leq \f{1}{2}\log(4d) + \log\log (4d) + \f{1}{2\sqrt{d}\log d} + 2 + \gamma
    \end{equation}
\end{corollary}

The corollary is a standard application of Riemann-Stieljes integrals. Let 

\[f(x) := \f{1}{2\Gamma(\nf{1}{2})\zeta(2)^{\nf{1}{2}}x}\paren{\f{1}{2}\log(4d) + \log\log (4d) + \f{1}{2\sqrt{d}\log d} + 2 + \gamma}^{\nf{1}{2}}.\]

We can compute the contribution of $d > D$ with another application of Riemann-Stieljes integration by parts applied to $f(x)dG(x)$. Note here that $d$
is the differential operator and $G(x)$, as in~\cref{thm:granchild}, counts the number of good $d$ with absolute value at most $x$.
We could likely improve upon~\cref{cor:pv-ineq} by using better bounds on $L(1,\chi)$~\cite{pintz},
but this is already sufficient to show~\crefrestated{cor:const-ub}, and in particular that $c_R > c_I$ under our hypothesis.  

\subsection{Putting everything together}

We apply the procedure of~\cref{subsubsec:good-d} to obtain the list of good~$d$ between $-10000$ and $50000$. Then, we apply~\cref{subsubsec:each-cd} with $Q := 1000$ to show that
\begin{equation}
    \label{eq:pet}
    \sum_{\substack{-10000 \leq d \leq 50000 \\ d \t{ good}}} \f{c_dc'_{d}}{|d|2^{\omega(d)}},
\end{equation}
where $c'_d$ is either $c'_{d,R}$ or $c'_{d,I}$. The lower bound on~\eqref{eq:pet} yields lower bounds on $c_R$ and $c_I$.

\constlb

We bound the contribution of good $d$ with $d$ outside $[-10000, 50000]$ using~\cref{assumption} as described in~\cref{subsubsec:tail}, and add it to the upper bound on~\eqref{eq:pet} to obtain upper bounds on $c_R$ and $c_I$. Furthermore, notice that~\cref{assumption} subsumes~\crefrestated{thm:red-ec} and thus we we no longer depend on the $abc$-conjecture. 

\begin{customcor}{E}
    \label{cor:const-ub:restated}
    Under~\cref{assumption} instead of the $abc$-conjecture,~\cref{thm:rix-const} holds with
    \[0.1255 \leq c_R \leq 0.1489 \quad\t{and}\quad 0.01109 \leq c_I \leq 0.03446.\]
    In particular, $c_R > c_I$ under this hypothesis. 
\end{customcor}

\bibliographystyle{plain}
\bibliography{egr-quad-ub}

\end{document}